\documentclass[11pt]{article}

\usepackage[top=1in,bottom=1.25in,left=1.25in,right=1.25in]{geometry}

\usepackage{graphicx}
\usepackage{ragged2e}
\usepackage{tabto}
\usepackage{amssymb}
\usepackage{float}
\usepackage{amsmath}
\usepackage{amsthm}
\usepackage{color}
\usepackage{url}
\usepackage{afterpage}
\usepackage[font=footnotesize]{caption}
\usepackage[font=footnotesize]{subcaption}
\usepackage{lmodern}
\usepackage[sf,bf,medium]{titlesec}

\usepackage{bm}
\usepackage{dsfont}
\usepackage{isomath}
\newcommand{\vct}{\vectorsym}

\newcommand{\lp}{\left(}
\newcommand{\rp}{\right)}

\usepackage[sort,compress]{cite}

\DeclareMathOperator\erf{Erf}

\newcommand\indic{\mathds 1}
\newcommand\bbR{\mathbb R}

\newcommand\bx{\vct{x}}

\newcommand\bc{\boldsymbol c}

\newcommand\bN{\boldsymbol N}

\newcommand\tfmm{T_{\textnormal{FMM}}}

\newtheorem{remark}{\sffamily Remark}
\newtheorem{definition}{\sffamily Definition}

\newcommand{\cO}{\mathcal O}

\newcommand{\cT}{\mathcal T}

\renewcommand{\phi}{\varphi}

\numberwithin{equation}{section}

\begin{document}


\begin{titlepage}

  \raggedleft
  {\sffamily Technical Report\\
    \today}
  
  \hrulefill

  \vspace{3\baselineskip}

  \raggedright
  {\huge \sffamily\bfseries A fast boundary integral method
    for high-order multiscale mesh generation}

  \vspace{3\baselineskip}
  {\sffamily \bfseries Felipe Vico}\footnote{Research supported in part by
  the Office of Naval Research under award numbers
  \#N00014-18-1-2307.}\\
\small \emph{Instituto de Telecomunicaciones y
  Aplicaciones Multimedia (ITEAM)\\
  Universidad Polit\`ecnica
  de Val\`encia\\
  Val\`encia, Spain 46022} \\
  \texttt{felipe.vico@gmail.com}

  \vspace{\baselineskip}
  \normalsize {\bfseries \sffamily Leslie Greengard}\footnote{Research supported in part by
  the Office of Naval Research under award numbers
  \#N00014-18-1-2307.}\\
  \small \emph{Courant Institute \tabto{2in} Center for Computational
    Mathematics\\
    New York University \tabto{2in} Flatiron Institute\\
    New York, NY, 10012 \tabto{2in} New York, NY 10010}\\
  \texttt{greengard@cims.nyu.edu}
  
  \vspace{\baselineskip}
  \normalsize {\bfseries \sffamily Michael O'Neil}\footnote{Research supported in part by
    the Office of Naval Research under award
    numbers~\#N00014-17-1-2059,~\#N00014-17-1-2451,
  and~\#N00014-18-1-2307.}\footnote{Corresponding author.}\\
\small \emph{Courant Institute\\
  New York University\\
    New York, NY 10012}\\
  \texttt{oneil@cims.nyu.edu}

  \vspace{\baselineskip}
  \normalsize {\bfseries \sffamily Manas Rachh}\\
  \small \emph{Center for Computational Mathematics\\
    Flatiron Institute\\ New York, NY 10010}\\
  \texttt{mrachh@flatironinstitute.org}

  \normalsize
  
\end{titlepage}

\begin{abstract}
  In this work we present an algorithm to construct an infinitely
  differentiable smooth surface from an input consisting of a
  (rectilinear) triangulation of a surface of arbitrary shape. The
  original surface can have non-trivial genus and multiscale
  features, and our algorithm has computational complexity which is
  linear in the number of input triangles. We use a smoothing kernel
  to define a function $\Phi$ whose level set defines the surface of
  interest.  Charts are subsequently generated as maps from the
  original user-specified triangles to $\bbR^3$. The degree of
  smoothness is controlled locally by the kernel to be commensurate
  with the fineness of the input triangulation.  The expression
  for~$\Phi$ can be transformed into a boundary integral, whose
  evaluation can be accelerated using a fast multipole method.  We
  demonstrate the effectiveness and cost of the algorithm with
  polyhedral and quadratic skeleton surfaces
  obtained from CAD and meshing software. \\


\end{abstract}

\tableofcontents

\newpage

\section{Introduction}

Over the past two decades, high-order accurate methods have been
developed for solving many of the partial differential equations
(PDEs) of classical mathematical physics in complicated geometries.
These include the equations that govern electromagnetics,
electrostatics, acoustics, fluid dynamics, and heat flow.  In two
dimensions, the state of the art is quite advanced, in the sense that
the combination of high-order representations of the domain,
high-order discretization methods, and fast algorithms have permitted
near machine precision accuracy using a modest number of unknowns.
For constant-coefficient, homogeneous equations, which can easily be
reformulated as boundary integral equations, it is sufficient to have
a high-order representation of the boundary itself.  After
discretization of the weakly-singular or singular integral operators
using high-order quadratures, see~\cite{rachh_2016} and the review
\cite{hao_2014}, the resulting linear systems can be solved in
optimal, or nearly optimal, time using fast multipole methods (FMMs)
or related fast algorithms (see, for example,
\cite{ggm1993,greengard-1987,helsing,martinsson-2005}).

The state of the art is quite different in the three-dimensional
setting. Standard boundary integral or boundary element methods often
use piecewise constant or linear approximations of boundary densities
defined on flat triangulated surfaces to discretize integral
operators.  Since fast algorithms
\cite{borm2003introduction,coifman1993fast,darve2004fast,greengard-1997,
  wideband3d,phillips1997,song1997multilevel,ying}
and high-order accurate quadrature schemes
\cite{wala2018qbx,Siegel2018,bremer-2015,bruno2001fast} are available,
it is the lack of suitable surface representations that has hampered
the development of accurate fast solvers in general geometries.  In
practice, the high-order geometries that can be accessed are those for
which an analytic parametrization is available (such as a deformed
sphere or torus) \cite{bremer-2015,oneil2018surface}, or ones for
which considerable effort has been expended in defining a smooth
surface.  Existing meshing algorithms which {\em do} generate
high-order curvilinear triangles or quadrilateral patches, such as
Gmsh~\cite{gmsh}, require a computer aided design-compatible geometry
as input.

In the present paper, we describe an algorithm to construct an
infinitely differentiable smooth surface from an input consisting of a
(rectilinear) triangulation of a surface of arbitrary shape.  The
original surface can have non-trivial genus and multiscale features,
and the running time of our algorithm is linear in the number of input
triangles.  We use a target-dependent (non-translation invariant)
smoothing  kernel whose integral
over the interior of the domain of interest defines a function $\Phi$
whose level set $\Phi = 1/2$ will be used as the definition of its
$C^\infty$ surface approximation.  Charts are then generated as maps
from the original user-specified triangles to $\bbR^3$.  The method is
able to handle highly non-uniform discretizations by ensuring that the
extent of smoothing is commensurate with the local fineness of the
triangulation.  By application of the divergence theorem,
the volume integral defining $\Phi$ can be converted into an integral over the
triangulated boundary itself; the resulting boundary integral can then
be evaluated using
the fast multipole method.  At the risk of stating the obvious, the
impetus behind our work is to reduce the overall cost of solving the
PDEs of mathematical physics.  High-order accuracy is critical,
for example, in large-scale wave propagation problems, but the ability
to achieve high-order accuracy is equally critical to enable automatic
adaptivity and robust error analysis even for non-oscillatory
problems.

There is an extensive literature on mesh generation,
mesh repair, surface parameterization, etc. Much of this has been driven
by the needs of computer graphics or visualization and we do not seek
to review the literature, but we do highlight earlier work that
is closest in spirit to the method presented below.
The problem of interest here involves taking as input
a water-tight flat triangulation, and creating an
infinitely smooth surface (and corresponding high-order triangulation)
which is a high-fidelity approximant of the low-order flat
triangulation. This problem clearly does not have a unique solution, 
and we seek simply to create {\em some} fixed surface which can serve as the 
underlying smooth surface to which 
a high-order accurate discretization is converging. 

Relevant earlier work includes~\cite{dapogny2014remesh}, in which
the authors propose a local method for mesh repair (adjusting the
size/aspect ratio of mesh
elements), based on local interpolation and adjustment of control
points. This algorithm is able to generate very high-quality meshes,
but is limited to curvilinear second-order triangles. Higher order
methods also exist that begin with 
flat triangulations, or even point clouds in~$\bbR^3$.  
One such scheme is ``moving least squares" surface reconstruction 
(see, for example,~\cite{fleishman2005}). This method
proceeds, as the name suggests, by computing a smooth surface based on a local
least-squares approximation of the data. It has the advantage of
being able to handle point clouds and noisy triangulations,
including those obtained from three dimensional scanning devices.
We should note that the algorithm below, in its current form, is less
general in that it does not handle noise. 
(We will return to this point in the concluding section.)

Perhaps most closely related to our scheme is a method
introduced in the 1990's, namely that of {\em convolution
surfaces}~\cite{bloomenthal1991,sherstyuk1999,sherstyuk1999design}.
Convolution surfaces were inspired by even earlier methods of surface
construction used in molecular modeling and animation, based on
computing iso-potential surfaces from point
sources~\cite{blinn1982}. The fundamental idea is to generate an implicit
surface by computing the convolution of a single radially symmetric Gaussian
with the characteristic function of the {\em surface}. In our scheme,
we define a level set function as the convolution of a Gaussian with the 
characteristic function of the {\em volume}, with a
variance carefully chosen in a target-specific fashion in order
to be commensurate with the fineness of the
nearby triangles. This enables a muti-scale form of smoothing. 
Another important class of methods is based on subdivision
surfaces~\cite{derose1998subdivision,zorin1996interpolating}, which
successively refine the triangulation as more and more smoothness is requested.
A scheme presented in~\cite{ying2004simple}
constructs an atlas based on large parameterized patches, taking
as input a polyhedral surface and using partitions of unity to blend them
together. Recently, a fast and robust tetrahedral meshing scheme
was developed~\cite{zorin2018}, 
which can also generate smooth surfaces, but its principal goal is the 
volume mesh, a topic we don't consider here.

Finally, we should note that important work has been done in
the area of isogeometric analysis
(IGA)~\cite{cottrell,hughes2005isogeometric,simpson2014acoustic} with
regard to boundary element and boundary integral equation methods. IGA
is a framework for dealing directly with CAD geometries, and provides
a robust set of tools for surface refinement, manipulation, and
discretization. There have been several PDE and integral equation
solvers constructed based on IGA, and significant progress has been
made toward coupling computer-aided design (CAD) systems with 
finite element methods. These schemes are complementary to the method 
presented below, which assumes no external information about the ``true"
underlying surface and takes as input only a flat triangulated mesh.

\begin{remark}
  A precursor of the algorithm presented here is the edge/corner
  rounding scheme of Epstein and O'Neil~\cite{epstein_2016}.  The
  basic observation in that work is the following: if corners of a
  polygon~$P \subset \bbR^2$
  are viewed locally as maps over the tangent line, then local convolution
  with a finite width ``bell" function results in a smooth curve which
  preserves convexity and leaves large flat regions unchanged. If the
  finite-width bell is replaced with a Gaussian which has decayed to
  magnitude~$\epsilon$, then the resulting geometry is
  an~$\epsilon$-accurate approximation to a true $C^\infty$ curve.
  While this scheme could be extended to an arbitrary
  three-dimensional geometry (the extension is, in fact, also
  described in~\cite{epstein_2016}), it would require the separate
  calculation of many local maps which would have to be patched
  together.  By reformulating the problem in terms of a volume
  integral with a $C^\infty$ kernel, no local maps need to be
  constructed, while retaining the useful features of the edge/corner
  rounding formalism.
\end{remark}

The paper is organized as follows: in Section~\ref{sec:smooth}, we construct
a multiscale mollifier for defining the smooth surface. 
In Section~\ref{sec:atlas}, we introduce
the mathematical foundations for defining an implicit surface as a level set 
and for constructing the corresponding atlas.
In Section~\ref{sec:algorithm}, we describe a fast algorithm for constructing
the level set and its charts.
Section~\ref{sec:examples} contains some numerical
examples demonstrating the efficiency and behavior of the algorithm,
and we conclude with a discussion of future work in
Section~\ref{sec:conclusions}.

\section{Smooth surfaces via multiscale mollification}
\label{sec:smooth}

Let $V \subset \bbR^3$ be a closed and bounded region with orientable
boundary~$S$, and let $K:\bbR^3 \to \bbR$ be a $C^\infty$ mollifier - that is, a
compactly supported function with
\begin{equation}
  \begin{aligned}
    \int_{\bbR^3} K(\bx) \, d\bx &= 1, \\
\lim_{\epsilon \rightarrow 0} \frac{K(\bx/\epsilon)}{\epsilon^3}  &= \delta(\bx),
\end{aligned}
\end{equation}
where $\delta(\bx)$ is the Dirac delta function. For
simplicity we will define~$K_\epsilon$ to be this scaled version
of~$K$:
\begin{equation}
K_\epsilon(\bx) = \frac{K(\bx/\epsilon)}{\epsilon^3} .
\end{equation}
It is well-known \cite{friedrichs} that
\begin{equation} \label{basemodel}
\begin{aligned}
\Phi(\bx) &= \int_{\bbR^3} K_\epsilon(\bx-\bx') \, \indic_V(\bx') \, d\bx' \\
 &= \int_{V} K_\epsilon(\bx-\bx')  \, d\bx'
\end{aligned}
\end{equation}
is itself a $C^\infty$ function, defined in all of~$\bbR^3$, where
$\indic_V$ is the characteristic function of $V$.  The level sets
$\Phi = C$ define smooth surfaces embedded in $\bbR^3$, and will be
denoted by $\Gamma_C$. The region bounded by~$\Gamma_C$ will be
denoted as~$\Omega_C$.  Moreover, interior regions of flat subsets
of~$S$ which are larger than the support of the kernel~$K_\epsilon$
are preserved along~$\Gamma_{1/2}$.  This is easy to show in the
one-dimensional case~\cite{epstein_2016}, and straightforward to prove
in higher dimensions. If~$V$ is convex, then~$\Gamma_{1/2}$ defines
a \emph{numerically-convex} region, in which the deviation from
convexity is on the scale of the square-root of the second moment of
the kernel (akin to the standard deviation).  Because of these
properties, in the remainder of this paper, we
will define the smooth surface via the level set~$\Gamma_{1/2}$.

Rather than using a mollifier in the strict mathematical sense
(i.e. one that is compactly supported), it is very
convenient for numerical purposes to replace~$K_\epsilon$ with 
the Gaussian kernel~$G_\sigma$:
\begin{equation}
  G_\sigma(\bx) = \frac{e^{-\lVert \bx \rVert^2/2\sigma^2}}{(2\pi
    \sigma^2)^{3/2}}, 
\end{equation}
and redefine~$\Phi$ using this kernel so that
\begin{equation} \label{basemodelg}
\Phi(\bx) = \int_{\bbR^3} G_\sigma(\bx-\bx') \, \indic_V(\bx') \, d\bx' 
 = \int_{V} G_\sigma(\bx-\bx')  \, d\bx' \, .
\end{equation}
Above,~$\Phi$ can be physically interpreted via heat flow: $\Phi(\bx)$
is the temperature at~$\bx \in \bbR^3$ at time~$\sigma^2/2$,
assuming that the temperature at
time zero is given by 1 for~$\bx \in V$ and zero otherwise.
Therefore, larger values of~$\sigma$ (i.e. time)
result in smoother distributions
of temperature.

The problem with this approach to surface generation as a
general-purpose tool, however, is that the single parameter $\epsilon$
or $\sigma$ determines a uniform length scale of smoothing.  In
geometries with multiscale features (involving triangles of vastly
different sizes), no single choice of $\epsilon$ or $\sigma$ can be
effective qualitatively or numerically:
the method would either fail to adequately smooth regions with large
triangles,  or it would wash out features in regions with small
triangles.  Thus, instead of~\eqref{basemodelg}, we will define~$\Phi$ by
\begin{equation} \label{basemodelk}
\Phi(\bx) 
 = \int_{V} G(\bx-\bx',\sigma(\bx))  \, d\bx' \, ,
\end{equation}
where
\begin{equation} \label{convkernel}
 G(\bx-\bx',\sigma(\bx)) = 
\frac{e^{-\lVert \bx-\bx' \rVert^2/2\sigma(\bx)^2}}{(2\pi \sigma(\bx)^2)^{3/2}},
\end{equation}
and $\sigma(\bx)$ is chosen to be commensurate with the size of the
triangle on~$S$ which is closest to~$\bx$.
We now turn to the actual construction of the non-constant~$\sigma$ above.

\subsection*{The multiscale mollifier}
\label{sec:sigma}

In order for the surface $\Gamma_{1/2}$ to be $C^{\infty}$, the
variance of the kernel~$G$, i.e. the function~\mbox{$\sigma = \sigma(\bx)$}
in~\eqref{convkernel}, must be $C^{\infty}$.
Furthermore,~$\sigma$ should be non-oscillatory in order to
preserve
the convexity of the domain as much as possible.  While there are many
possibilities of constructing such a function,  we use
the following formulation, which works well in
practice for surfaces with multiscale features.

Suppose that~$S$, which we will also refer to as the \emph{skeleton
  surface}, is given as the union of $M$ triangles $T^{j}$,
$j=1,2,\ldots M$. Let $\bc_{j}$ denote the centroid of $T^{j}$ and let
$\sigma_{j}= D_{j}/\lambda$, where~$D_{j}$ is the diameter of $T^{j}$
(i.e. the diameter of the smallest ball which encloses~$T^j$), and
$\lambda$ is a free parameter.
 We then  define~$\sigma$ as:
\begin{equation}\label{sigma2}
\sigma(\bx)=\frac{\sum_{j=1}^{M} \sigma_j  \, 
e^{-{\lVert\bx-\bc_j\rVert^2} / {2\sigma_0^2}}}
{\sum_{j=1}^{M} e^{-\lVert\bx-\bc_j\rVert^2 / 2\sigma_0^2}},
\end{equation}
where $\sigma_0$ is a second free parameter.
Note that the dominant term in the sum in the numerator is due to the
nearest centroid~$\bc_j$ to the point~$\bx$; other contributions
decay exponentially fast.
It is useful to set these parameters such that the influence on a
triangle $T$ due to triangles $T^j$ of the same size, but outside $T$'s
immediate nearest neighbors, is nearly negligible. 
Ultimately,~$\lambda$ controls the smoothness of the
final surface and~$\sigma_0$ controls the smoothness of~$\sigma$
itself. In practice, reasonable choices for the free parameters 
  above are
  \begin{equation}
    \lambda \approx 2.5, \qquad \sigma_0 \approx \sqrt{5} \max_j D_j.
  \end{equation}
Using these values, the influence on~$\sigma$ of triangles separated by a
diameter is roughly~$0.5$ (due to the choice of~$\sigma_0$), and
$\sigma$ has decayed to approximately~$10^{-9}$ at a distance
of~$2.5$ diameters from the centroid~(as determined
by~$\lambda$). Such a choice enables a straightforward splitting of
the near field and far field of triangle~$T$ when embedded in an
octree data structure (see Section~\ref{sec:fmm}).
With the above observations in mind, and the formula for
computing~$\Phi$ in~\eqref{basemodelk}, we now turn to the
construction of an atlas defining the~$C^\infty$ as a collection of
charts from the input skeleton triangulation.

\section{Constructing an atlas}
\label{sec:atlas}

In this section, we present a method for constructing an atlas (i.e. a
collection of charts or parameterizations) for $\Gamma_{1/2}$.  Recall
that the skeleton surface $S$ is given as the union of~$M$
triangles~$T^j$, $j=1,\ldots,M$.  For each triangle $T^j$, we denote
its three vertices by $\{ \vct{P}^j_1,\vct{P}^j_2,\vct{P}^j_3 \}$.  We
assume this triangulation of $S$ is water-tight, positively oriented, and
conforming (by conforming, we mean here that triangles which touch
either share a common edge or meet only at a triangle vertex.  It is
straightforward to extend our construction to other configurations,
including quadrilateral patches, non-conforming triangulations, etc.).

Each skeleton triangle can be parameterized as:
\begin{equation}\label{eq:tparam}
  \vct{T}^j(u,v) = \vct{P}^j_1 + u \lp \vct{P}^j_2 - \vct{P}^j_1 \rp
  + v \lp \vct{P}^j_3 - \vct{P}^j_1 \rp,
\end{equation}
with a local basis for~$\mathbb R^3$ defined by
\begin{equation}
\vct{T}^j_u = \vct{P}^j_2 - \vct{P}^j_1, \qquad
\vct{T}^j_v = \vct{P}^j_3 - \vct{P}^j_1, \qquad
\vct{N}^j =  \vct{T}^j_u \times \vct{T}^j_v.
\end{equation}
In \eqref{eq:tparam}, $(u,v)$ must lie in the standard
simplex triangle  
\begin{equation} \label{stdsimplex}
T_0 = \{ (u,v) | \  u \geq 0, v \geq 0, (u+v) \leq 1 \}.  
\end{equation}
Note that above, surfaces in $\bbR^3$, such as $T^j$, are denoted
using standard weight font, and vectors, such as~$\vct{T}^j$, are
denoted in boldface.  The normalized versions of the above vectors
will be given
as~$\hat{\vct{T}}^j_u$,~$\hat{\vct{T}}^j_v$,~$\hat{\vct{N}}^j$
(in general, the vectors~$\hat{\vct{T}}^j_u$
and~$\hat{\vct{T}}^j_v$ are not orthogonal).

We now wish to construct a {\em mapping direction}~$\vct{H}^j$ along
triangle $T^j$ such that the level surface~$\Gamma_{1/2}$ can be
parametrized as the union of curved triangular patches $\Gamma^j$, and
therefore~$\Gamma_{1/2} = \cup_{j=1}^{M} \Gamma^{j}$, through the
charts
\begin{equation}\label{eq:h}
  \vct{x}^j(u,v) = \vct{T}^j(u,v) + h^j(u,v) \, \vct{H}^j(u,v),
\end{equation}
where $\vct{x}^{j} : T_{0} \to \Gamma^{j}$, and $h_{j}(u,v)$ is
determined so that $\vct{x}^{j} \subset \Gamma_{1/2}$.  For this to be
possible, the mapping direction must not only be continuous along~$S$,
but imply a bijective mapping from~$S$ to $\Gamma_{1/2}$.  Thus, we
cannot use the triangle normals $\vct{N}^j$ themselves, since they are
discontinuous across triangle edges and would yield gaps in the
atlas defining the smooth surface. However, a continuous mapping
direction, which we refer to as the {\em pseudonormal} vector field
can be constructed on each triangle as follows.

\begin{definition}
  Let $\vct{P}$ be a vertex on the skeleton surface $S$.  The set of
  triangles for which~$\vct{P}$ is a common vertex is denoted
  by~$\cT(\vct{P})$.  The {\em vertex angle} $\theta_{\vct{P}}(T)$ of
  triangle~$T \in \cT(\vct{P})$ is defined to be the interior angle
  of~$T$ at this vertex.
\end{definition}

\begin{figure}[t]
  \centering
  \includegraphics[width=.8\linewidth]{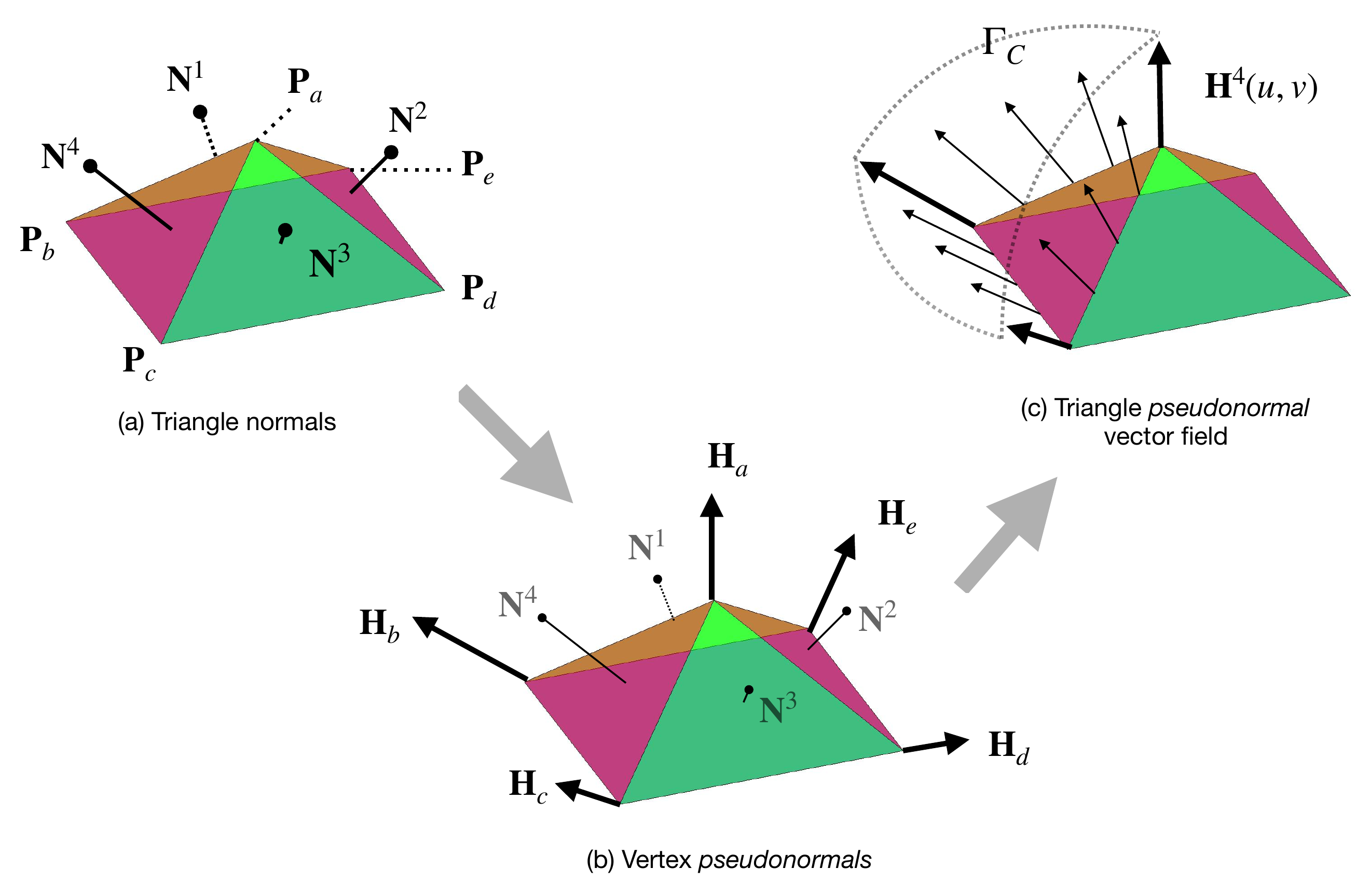}%
  \caption{In (a), we show vertex $\vct{P}_a$ and the other vertices
    $\vct{P}_b, \vct{P}_c, \vct{P}_d,$ and $\vct{P}_e$ of the four
    triangles impinging on it. The normal vectors to the four
    triangles are denoted by $\vct{N}_1, \dots, \vct{N}_4$. In (b), we
    show the vertex pseudonormals, constructed according to
    \eqref{mkvertexpn}. In (c), we plot the pseudonormal vector field
    $\vct{H}^4(u,v)$ along the edges of triangle 4, constructed
    according to~\eqref{mkpseudovec}. We also show a portion of the
    level set $\Gamma_{1/2}$, which will be represented
    using~\eqref{eq:h}.  (\emph{Note: The triangles depicted here are
      not a good approximation of an underlying smooth surface, and
      serve only as an illustration.})}
\label{fig:normalverts}
\end{figure}

\begin{definition}
Let $\vct{P}$ be a vertex on the skeleton surface $S$,
and let $A(\vct{P})$ denote the total vertex angle
at $\vct{P}$:
\begin{equation}
A(\vct{P}) =  \sum_{T \in \cT(\vct{P})} \theta_{\vct{P}}(T).
\end{equation}
The {\em vertex pseudonormal} $\vct{H}(\vct{P})$ is then given by
\begin{equation}  \label{mkvertexpn}
  \vct{H}(\vct{P}) = \frac{1}{ A(\vct{P})} \, 
  \sum_{T \in \cT(\vct{P})} 
  \theta_{\vct{P}}(T) \,  \hat{\vct{N}}(T),
\end{equation}
where~$\hat{\vct{N}}(T)$ is the unit normal along triangle~$T$.
Furthermore, on triangle~$T^{j}$, let the vertex psuedonormals
corresponding to its
vertices~$\{ \vct{P}^j_1,\vct{P}^j_2,\vct{P}^j_3 \}$ be denoted
by~$\{ \vct{H}^j_1,\vct{H}^j_2,\vct{H}^j_3 \}$.  Then, the
\emph{pseudonormal vector field} for triangle $T^j$ is defined by the
convex combination:
\begin{equation} \label{mkpseudovec}
  \vct{H}^j(u,v) = \vct{H}^j_1 + u \lp \vct{H}^j_2 - \vct{H}^j_1 \rp
  + v \lp \vct{H}^j_3 - \vct{H}^j_1 \rp,
\end{equation}
where $(u,v)$ lies on the standard simplex triangle.
See Figure~\ref{fig:normalverts} for a depiction.
We will refer to the function $h^j$ 
in~\eqref{eq:h} as the {\em pseudonormal distance function}. 
\end{definition}

In short, the vertex pseudonormal is a weighted average of the normals
of all triangles impinging on that particular vertex, and the
pseudonormal vector field is continuous on $S$.  We note that there
are many possible definitions for the weighted average used to define
the vertex pseudonormals.  For any such choice, the construction
in~\eqref{mkpseudovec} yields a continuous vector field on~$S$.

Turning now to the pseudonormal distance function~$h^j$, it is clearly
determined point-wise on triangle $T^j$ as the function~$h^j$ which
satisfies the scalar equation
\begin{equation}\label{eq:newton}
  \Phi\lp 
  \vct{T}^j(u,v) + h^j(u,v) \, \vct{H}^j(u,v) \rp - \frac{1}{2} = 0.
\end{equation}
The solution to~\eqref{eq:newton} can easily be found using Newton's
method, for which we will require the evaluation of the directional
derivative of~$\Phi$:
\begin{equation}\label{eq:newtonderiv}
  \frac{\partial \Phi}{\partial h^j}(u,v) = 
  \vct{H}^j(u,v) \cdot \nabla  \Phi(\vct{T}^j(u,v) + h^j(u,v) \, \vct{H}^j(u,v)). 
\end{equation}
An expression for~$\nabla \Phi$ can be determined directly from its
integral representation in~\eqref{basemodelk}.
Once $h^j$ is known, 
the final surface is obtained as the
collection of the corresponding charts~$\vct{x}^j$, $j=1,2,\ldots M$.  
The previous discussion implies that the function~$h^j$ has
been computed for every possible~$u,v$, or that it is computed when needed
for each~$u,v$. Shortly, in Section~\ref{subsec:hoapprox}, we
detail exactly where and how~$h^j$ is evaluated.

\subsection{Local coordinates}
\label{sec:local}

Our motivating discussion in the introduction laid out the
requirements for high-order boundary integral equation solvers in
three dimensions. If the surface constructed via the method of this
paper is to be used as a high-order geometry processing algorithm, we
will also require suitable expressions for the area element and surface tangents
on each patch~$\Gamma^j$ of the level surface $\Gamma_{1/2}$.
Recall that the chart for the patch~$\Gamma^j$ is given by
\begin{equation}\label{chart}
    \vct{x}^j(u,v) = \vct{T}^j(u,v) + h^j(u,v) \, \vct{H}^j(u,v).
\end{equation}
Taking the partial derivative with respect to~$u$ above, we have
\begin{equation}\label{u_vector}
  \begin{aligned}
    \vct{x}^j_u(u,v) &=  \frac{\partial \vct{x}^j}{\partial u}(u,v) \\
    &= \vct{T}^j_u(u,v) + h^j(u,v) \,
  \vct{H}^j_u(u,v)   + \frac{\partial h^j}{\partial
    u}(u,v)
  \, \vct{H}^j(u,v).
  \end{aligned}
\end{equation}
Here,~$\partial h^j / \partial u$ is unknown, since it was only~$h^j$
that was determined via Newton's method.  However, using the fact
that~$\vct{x}^j_u$ is tangent to the
surface~$\Gamma_{1/2}$, and that~$\nabla \Phi$ is normal to~$\Gamma$
(by the definition of a level set), we have
\begin{equation}\label{U_vector_2}
  \frac{\partial h^j}{\partial u} = - \frac{\lp \vct{T}^j_u + h^j \, \vct{H}^j_u
  \rp \cdot \nabla \Phi }{\vct{H}^j \cdot \nabla \Phi}.
\end{equation}
A similar calculation can be used to obtain
\begin{equation}\label{v_vector}
    \vct{x}^j_v
    = \vct{T}^j_v + h^j \,
  \vct{H}^j_v   + \frac{\partial h^j}{\partial
    v}
  \, \vct{H}^j
\end{equation}
where
\begin{equation}\label{v_vector_2}
  \frac{\partial h^j}{\partial v} = - \frac{\lp \vct{T}^j_v + h^j \,
    \vct{H}^j_v
  \rp \cdot \nabla \Phi }{\vct{H}^j \cdot \nabla \Phi}.
\end{equation}
The area element on~$\Gamma^j$ is easily computed from the definition
\begin{equation}\label{eq:da}
 da = \lvert \bx^j_u \times \bx^j_v \rvert \, du \, dv, 
\end{equation}
and the normal vector, as noted above, is given by~$\nabla \Phi$.

Depending on the application~\cite{oneil2018surface},
it may be necessary to compute higher
order derivatives of the surface parameterization. These derivatives
can be obtained via a computation similar to that used in obtaining
first derivatives.
Consider the calculation of the second order partial
derivative $\vct{x}^j_{uu}(u,v)$.  Taking the derivative
of both sides of~\eqref{u_vector} with respect to~$u$ we have:
\begin{equation}\label{uu_vector}
  \begin{aligned}
    \vct{x}^j_{uu}(u,v) &= \vct{T}^j_{uu}(u,v) + h^j(u,v) \,
  \vct{H}^j_{uu}(u,v)   + 2 h_u^j(u,v)+ h_{uu}^j(u,v)
  \, \vct{H}^j_u(u,v).
  \end{aligned}
\end{equation}
In the expression above
the only unknown term is~$h_{uu}^j(u,v)$. This term can be obtained
by differentiating both sides of equation~\eqref{eq:newton} twice:
\begin{multline}\label{uu_vector_aux1}
  \lp \vct{T}^j_u+h^j_u\vct{H}^j +h^j\vct{H}^j_u \rp^T \nabla\nabla\Phi
  \lp\vct{T}^j_u+h^j_u\vct{H}^j +h^j\vct{H}^j_u \rp + \\
  +\nabla\Phi\cdot(\vct{T}^j_{uu}+2h^j_u\vct{H}^j_u + h\vct{H}^j_{uu}
  +h_{uu}\vct{H}^j)=0,
\end{multline}
where~$\nabla\nabla \Phi$ is the Hessian matrix of the
function~$\Phi$. Then, solving for~$h^j_{uu}$ we get:
\begin{multline}\label{uu_vector_aux2}
  h^j_{uu}=-\frac{\lp \vct{T}^j_u+h^j_u\vct{H}^j +h^j\vct{H}^j_u
    \rp^T
    \nabla\nabla\Phi
  \lp\vct{T}^j_u+h^j_u\vct{H}^j +h^j\vct{H}^j_u \rp} {\nabla\Phi\cdot\vct{H}^j}\\
-\frac{\nabla\Phi\cdot(\vct{T}^j_{uu}+2h^j_u\vct{H}^j_u
  +h\vct{H}^j_{uu})}{\nabla\Phi\cdot\vct{H}^j}.
\end{multline}
Similar expressions for~$\vct{x}^j_{uv},\vct{x}^j_{vv}$, or any
arbitrary higher order derivative of the surface, can be derived (of
course the algebra becomes rather unwieldy rather quickly). Using
these higher order derivatives, it becomes possible to obtain the
second fundamental form of the surface, Gaussian/mean curvatures,
Christoffel symbols, and other quantities of interest in
differential geometry. These quantities are often used in time
dependent PDEs to determine the evolution of free boundaries.
Special care may need to be used when
numerically computing these higher derivatives, but schemes are
necessarily case-dependent.

\subsection{High-order approximation}
\label{subsec:hoapprox}

We now seek to develop a high-order approximation of the
surface~$\Gamma_{1/2}$ as a collection of piecewise smooth charts.  To
this end, we will approximate each component of each
chart~$\vct{x}^j : {T}_{0} \to \Gamma^{j}$ in~\eqref{eq:h} as a
polynomial in $u,v \in {T}_0$ of total degree $p$. Each chart will be
referred to as a $p$th-order curvilinear triangle, or $p$th-order
approximant.  This can be achieved by a method analogous to polynomial
interpolation/approximation in one dimension: the
function~$\vct{x}^{j}$ will be sampled at the $p$th-order
Vioreanu-Rokhlin nodes on ${T}_{0}$, which provide stable
interpolation formulae for high-order polynomial approximation on the
simplex~${T}_{0}$~\cite{vioreanu_2014}. Polynomial approximation and
interpolation is performed using an orthogonal basis of Koornwinder
polynomial~\cite{koornwinder_1975}, analogous to using Legendre
polynomials on the interval~$[-1,1]$. Note that in contrast to some
computational geometry procedures for determining surfaces, no effort
is made to match values or derivatives at interfaces, but rather only
to approximate the function to high accuracy on each panel. (And
therefore any gaps that may appear in our approximation
of~$\Gamma_{1/2}$ can be controlled and made to be as small as
desired.)

To this end, we first recall that there are~$n_p = (p+1)(p+2)/2$
polynomials of two variables with total degree~$\leq p$.
Let $(u_{i},v_{i})$,
$i=1,2,\ldots n_{p}$ denote the $p$th-order
Vioreanu-Rokhlin nodes.
As detailed earlier, in order to evaluate~$\vct{x}^{j}(u_{i},v_{i})$
it is necessary to compute the corresponding pseudonormal 
distance~$h^{j}(u_{i},v_{i})$ which satisfies:
\begin{equation}
  \Phi(\bx^j (u_{i},v_{i})) = \Phi\lp \vct{T}^j(u_i,v_i)
  + h^j(u_i,v_i) \, \vct{H}^j(u_i,v_i) \rp = \frac{1}{2}.
\end{equation}
To simplify the notation, let~$h_{ji} = h^j(u_i,v_i)$; then,
letting~$h^{(k)}_{ji}$ be the~$k$th iterate for
computing $h_{ji}$ in Newton's method, we have
\begin{equation}\label{newtoneq}
    h_{ji}^{(k+1)} = h_{ji}^{(k)} -
\frac{\Phi\lp \vct{T}^j(u_i,v_i) + h_{ji}^{(k)} \, \vct{H}^j(u_i,v_i) \rp - 1/2}
{\vct{H}^j(u_{i},v_{i}) \cdot \nabla  \Phi \lp \vct{T}^j(u_{i},v_{i})
  + h_{ji}^{(k)}\, \vct{H}^j(u_{i},v_{i}) \rp},
\end{equation}
where we make use of~\eqref{eq:newtonderiv} for
$\partial{\Phi}/\partial{h}$. We initialize the Newton iteration with
$h_{ji}^{(0)} = 0$.  Once Newton's method has converged, the local
coordinate system and metric tensor along the surface, with respect to
the local parameterizations on~$T^j$, can be calculated using
the expressions in Section~\ref{sec:local}.

\begin{figure}[t]
  \centering
  \includegraphics[width=0.4\linewidth]{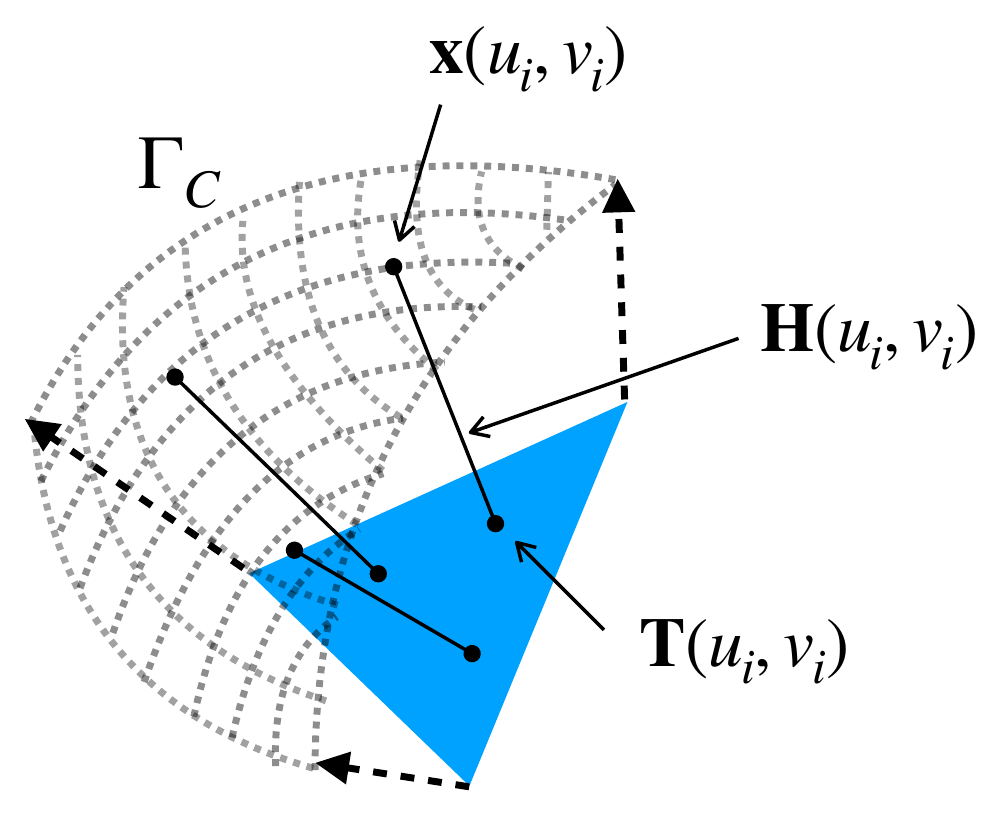}
  \caption{The pseudonormal direction vector $\vct{H}(u_i,v_i)$
    emanating from the quadrature node $\vct{T}(u_i,v_i)$ intersects
    the level surface $\Gamma_C$ at $\vct{x}(u_i,v_i)$. The distance
    along $\vct{H}(u_i,v_i)$, denoted by $h_{ji}$ is determined using
    Newton's method.  }
\label{fig:newton} 
\end{figure}

\section{Computing the level set function}
\label{sec:algorithm}

Each Newton step for evaluating~$h^{(k)}_{ji}$, used in eventually
computing the point~\mbox{$\bx^{j}(u_{i},v_{i}) \in \Gamma^j$}, requires the
evaluation of~$\Phi$ and its gradient at $n_p M$ target locations
($n_p$ nodes on each of~$M$ patches).  Thus, in order for the method
to be practical, we need to be able to compute the volume integral in
\eqref{basemodelk} accurately and rapidly, preferably with a
computational cost scaling as~$\cO \lp n_p M \rp$, i.e. the number of
interpolation points on~$\Gamma_{1/2}$.
 One option would
be to design a variant of the fast Gauss transform (FGT) that combines
the hierarchical data structure of \cite{wang2018} with the variable
scale FGT of \cite{strain91}. This, however, would require a
surface-conforming volume mesh for the region~$V$ and a rather
complicated implementation of the fast algorithm.  Instead, we will
recast~\eqref{basemodelk} as a boundary integral and discuss fast
algorithms for evaluating the reformulated version of~$\Phi$.

\subsection{Surface integral representation}
\label{sec:surface}

The integral representation in~\eqref{basemodelk} can be reformulated
as a boundary integral along~$S$ by a straightforward application of
the divergence theorem. First, we observe that
\begin{equation}
  \Delta_{\bx'} \lp 
  \frac{\erf \lp \frac{\lVert\bx-\bx'\rVert}{\sqrt{2}\sigma(\bx)} \rp}{4\pi \lVert\bx-\bx'\rVert} 
  \rp = - \frac{e^{-\lVert\bx-\bx'\rVert^2/2\sigma(\bx)^2}}{\lp 2\pi \sigma(\bx)^2 \rp^{3/2}}, 
\end{equation}
where
\begin{equation}
\erf(r) = \frac{2}{\sqrt{\pi}} \int^r_0 e^{-t^2} \, dt.
\end{equation}
From the preceding identity, we can rewrite~$\Phi$ from \eqref{basemodelk} as
\begin{equation}\label{eq:surf}
  \begin{aligned}
    \Phi(\vct{x}) &= \int_V 
\frac{e^{-\lVert\vct{x}-\vct{x}'\rVert^2/2\sigma(\bx)^2}}{\lp 2\pi \sigma(\bx)^2 \rp^{3/2}}
    \, d\vct{x}' \\
&= - \int_V \Delta_{\bx'} \lp \frac{\erf \lp \frac{\lVert\vct{x}-\vct{x}'\rVert}{\sqrt{2}\sigma(\bx)}\rp}{4\pi \lVert\vct{x}-\vct{x}'\rVert}
        \rp \, d\vct{x}' \\
    &= - \int_S \hat{\vct{N}}(\bx') \cdot \nabla_{\bx'}
     \lp \frac{\erf \lp \frac{\lVert\vct{x} - \vct{x}'\rVert}{\sqrt{2}\sigma(\bx)} \rp}
{4\pi \lVert\vct{x}-\vct{x}'\rVert}
    \rp \, da(\vct{x}') \\
    &= \int_S \psi_\sigma(\bx,\bx') \, da(\bx') \, ,
  \end{aligned}
\end{equation}
where
\begin{equation}\label{eq:kerpsi}
  \begin{aligned}
    \psi_\sigma(\bx,\bx') &= 
      -\hat{\vct{N}}(\vct{x}') \cdot \lp \vct{x} -\vct{x}' \rp \lp
    \frac{\erf \lp \frac{\lVert\vct{x} - \vct{x}'\rVert}{\sqrt{2}\sigma(\bx)} \rp
    }{4\pi\lVert\bx-\bx'\rVert^3}
    -\frac{\sqrt{\frac{2}{\pi}}e^{-\lVert\bx-\bx'\rVert^2/2\sigma(\bx)^2}}{4\pi\sigma(\bx)\lVert\bx-\bx'\rVert^2}
    \rp, \\
    da(\bx') &= \left\lVert \vct{T}_u \times \vct{T}_v \right\rVert \, du \,dv,
  \end{aligned}
\end{equation}
and~$\hat{\vct{N}}$ denotes the unit outward normal along the skeleton
surface. The index~$j$ denoting individual triangles on the skeleton
mesh has been suppressed in the above expressions.
Note that we have made use of the divergence theorem 
to obtain the third line of~\eqref{eq:surf}.

In order to compute~$\nabla \Phi$, we will also need to 
evaluate~$\nabla \psi_\sigma$.  From~\eqref{eq:surf}
and~\eqref{eq:kerpsi}, we have
\begin{multline}\label{gradPhi}
  \int_S \nabla_{\bx} \psi_\sigma(\bx,\bx') \, dS(\bx') \\
    =\int_S-\hat{\bN}(\bx') \lp
    \frac{\erf \lp \frac{\lVert\bx-\bx'\rVert}{\sqrt{2}\sigma(\bx)} \rp
    }{4\pi\lVert\bx-\bx'\rVert^3}
    -\frac{\sqrt{\frac{2}{\pi}}e^{-\frac{\lVert\bx-\bx'\rVert^2}{2\sigma^2(\bx)}}}
    {4\pi\sigma(\bx)\lVert\bx-\bx'\rVert^2}    \rp \, dS(\bx')\\
    +\int_S (\bx-\bx') \, F(\bx,\bx',\sigma)
    \lp (\bx-\bx')\cdot \hat{\bN}(\bx') \rp \,  dS(\bx') \\
    -\int_S
    \frac{e^{-\frac{\lVert\bx-\bx'\rVert^2}{2\sigma^2(\bx)}}}{\sigma^4(\bx)\pi^{3/2}2\sqrt{2}}
    \lp (\bx-\bx') \cdot \hat{\bN}(\bx') \rp \, 
    \nabla_{\bx}\sigma(\bx) \, dS(\bx'),
\end{multline}
where 
\begin{equation}
F(\bx,\bx',\sigma) = 
\frac{e^{-\frac{r^2}{2\sigma^2(\bx)}}\sqrt{2}(r^3+3
  r\sigma^2(\bx))-3\sqrt{\pi}\sigma^3(\bx)\erf{ \lp \frac{r}{\sqrt{2}
      \sigma(\bx)} \rp} }{4 \pi^{3/2} r^5\sigma(\bx)^3}
\end{equation}
and $r = \lVert\bx - \bx'\rVert$.

\subsection{Numerical evaluation of $\sigma$}
\label{sec:evalsig}

Given the parameters~$\lambda$ and~$\sigma_0$ which determine the shape
of~$\sigma$ in~\eqref{sigma2}, evaluation is straightforward.
Having sorted the skeleton triangles into an octree data structure
based on their centroids, the evaluation of the function~$\sigma$ is a
purely local calculation due to the decay of the Gaussian
kernels. Therefore, no fast algorithm is needed for its evaluation,
merely maintaining an octree data structure is sufficient.

\subsection{Numerical evaluation of $\Phi$}
We turn now to the numerical evaluation of the level set function
$\Phi(\bx)$ and its gradient using~\eqref{eq:surf}
and~\eqref{gradPhi}.  Since the kernels defining the integrands in
\eqref{eq:surf} and \eqref{gradPhi} are smooth, we require only a
high-order quadrature rule for smooth functions on triangles.  For this
purpose, we will make use of Vioreanu-Rokhlin
rules~\cite{vioreanu_2014}, which also serve as interpolation nodes for
polynomials on triangles.
\begin{remark}
  The Vioreanu-Rokhlin (VR) quadrature rules are \emph{Gaussian-like} in that they
  integrate more functions than there are nodes in the quadrature. A
  perfect $n$-point Gaussian rule in two variables would exactly
  integrate $3n$ functions (as there are $3n$ parameters).  In our
  case, let $n_q = (q+1)(q+2)/2$ denote the number of polynomials of
  total degree less than or equal to~$q$ in two variables. Using~$n_q$
  nodes and weights, the Gaussian-like VR rules integrate
  exactly polynomials of $u,v$ of total degree less than or equal to
  $q'$ for some $q < q'$, and where $n_{q'} < 3n_q$.  The relationship
  between $q$ and $q'$ is somewhat complicated, and we refer the
  reader to the original paper for details.  For illustration, we note
  from~\cite{vioreanu_2014}, that with $q=4$, there are 15
  interpolation/quadrature nodes on the simplex, and that the
  VR quadrature rules integrate all polynomials of total
  degree less than or equal to $q' = 7$, of which there are 36 such
  functions.
\end{remark}

Once each skeleton triangle has been discretized using an
$n_q$-point quadrature rule with nodes~$\{ u_k,v_k\}$ 
and weights~$\{ w_k \}$ on~${T}_0$, we replace the 
surface integral representation for~$\Phi$ in~\eqref{eq:surf}
with the discrete sum
\begin{equation} \label{discrete_int}
  \Phi(\bx) \approx \sum_{j = 1}^M \sum_{k = 1}^{n_q} w_k \, \psi_{\sigma} \lp
  \bx-\vct{T}^j(u_k,v_k)\rp,
\end{equation}
where~$\vct{T}^j$ is given in~\eqref{eq:tparam}.
This approximation is defined in all of~$\bbR^3$, and its 
level set~$\Gamma_{1/2}$ (an approximation to the \emph{true} level
set) defines an analytic smooth surface.

\begin{remark}
  Note that the $n_{q}$ nodes (quadrature nodes) for computing the
  integral in~\eqref{discrete_int} are separate from the $n_{p}$ nodes
  (discretization nodes) for obtaining a $p$th-order representation of
  the charts in~\ref{subsec:hoapprox}.  The quadrature nodes
  determine the continuous level set $\Gamma_{1/2}$ and the
  discretization determines the order of accuracy of the
  representation of $\Gamma_{1/2}$. In practice, it is often the case
  that~$q$ is chosen to be large so as to ensure the accurate
  evaluation of the integral, and~$p$ is of modest size.
\end{remark}

\subsection{Fast multipole acceleration}
\label{sec:fmm}

Referring to~\eqref{newtoneq}, at the $k$th Newton iterate for computing
the pseudonormal distance $h_{ji}$, we need to evaluate
\begin{equation}
  \begin{aligned}
    \Phi \lp \bx_{ji}^{(k)} \rp &=  \sum_{\ell=1}^{M} \sum_{m=1}^{n_{q}} w_{m}
    \, 
    \psi_{\sigma} \lp \bx_{ji}^{(k)} - \vct{T}^{\ell}(u_{m},v_{m}) \rp,\\
    \nabla \Phi  \lp \bx_{ji}^{(k)} \rp &= \sum_{\ell=1}^{M} \sum_{m=1}^{n_{q}}
    w_{m} \, \nabla \psi_{\sigma} \lp \bx_{ji}^{(k)} -
    \vct{T}^{\ell}(u_{m},v_{m}) \rp ,
  \end{aligned}
\end{equation}
where 
\begin{equation}
  \bx_{ji}^{(k)} = \vct{T}^{j}(u_{i},v_{i}) + h_{ji}^{(k)} \,
  \vct{H}^{j}(u_{i},v_{i}).
\end{equation}
Thus, computing~$\Phi(\bx_{ji}^{(k)})$ and~$\nabla \Phi (\bx_{ji}^{(k)})$
boils down to computing $N$-body sums involving~$n_{q}M$
sources and~$n_{p}M$ targets for the kernels given by~$\psi_{\sigma}$
and~$\nabla \psi_{\sigma}$.  This computation can be
accelerated using a variety of fast algorithms and evaluated
in~$\cO \lp (n_p + n_{q})M \rp$ CPU time.  In this section, we
briefly describe one such approach, based on the fast multipole method
for Laplace's equation.  We refer the reader
to~\cite{fmm2,greengard-1997} for a thorough description of that method
and to a sampling of the literature
\cite{borm2003introduction,darvebb,wideband3d,phillips1997,song1997multilevel,ying}
for related fast summation schemes. The exact FMM code used in our
algorithm is based on the FMM3D library at
\url{github.com/flatironinstitute/fmm3d}.

We first observe that once~$r = \lVert \bx-\bx' \rVert \geq 8 \sigma$,
\begin{equation}
  \begin{aligned}
    \lvert \erf(r/(\sqrt{2} \sigma)) - 1 \rvert &\leq 10^{-14}, \\
    e^{-r^2/(2\sigma^2)} &\leq 1.5 \cdot 10^{-14}.
  \end{aligned}
  \end{equation}
Thus, at distances greater than~$8\sigma$, 
the kernel~$\psi_\sigma$ can be approximated to near double
precision as
\begin{equation} \label{kerapprox}
  \psi_\sigma(\bx-\bx') \approx -\hat{\vct{N}}(\bx') \cdot \nabla_{\bx'} \left(
    \frac{1}{4\pi \lVert\bx-\bx'\rVert }\right).
\end{equation}
The above expression is simply the kernel of the double layer
potential for the Laplace operator, to which FMMs for Laplace
potentials apply directly.  More precisely, for any target point
$\bx \in \bbR^3$, we may write
\begin{equation} \label{kerdecomp}
  \Phi(\bx) \approx -\int_{S\setminus B_R(\bx)} \hat{\vct{N}}(\bx')
  \cdot \nabla_{\bx} \lp
    \frac{1}{4\pi \lVert\bx-\bx'\rVert } \rp \, da(\bx') +
  \int_{S \cap B_R(\bx)} \psi_\sigma(\bx,\bx') \, da(\bx'),
\end{equation}
where~$B_R(\bx)$ is a ball of radius~$R$ centered at~$\bx$ with
radius~$R = R\lp \sigma(\bx) \rp$ chosen so that~\eqref{kerapprox} is
correct to the desired precision. The first term can be computed for
all~\mbox{$\bx \in \Gamma_{1/2}^{(k)}$} by the FMM in linear time. Here,
$\Gamma_{1/2}^{(k)}$ denotes the $k$th approximation to the
surface~$\Gamma_{1/2}$ obtained during the Newton iteration described
in Section~\ref{subsec:hoapprox}.  The second term
in~\eqref{kerdecomp} is a purely local calculation,
and can be carried out directly without
the need for a fast algorithm.  Assuming that $\sigma(\bx)$ is
approximately of the same length scale as the nearest triangle and
that the triangulation itself is multiscale but suitably graded, it
is straightforward to show that the total cost for evaluating all
such local
interactions is $\cO \lp M(n_p+n_q) \rp$.

A naive implementation for evaluating $\psi_{\sigma}(\bx,\bx')$,
however, is subject to catastrophic cancellation when $\bx$ is close
to $\bx'$; since the expression for $\psi_{\sigma}$ in
\eqref{eq:kerpsi} involves the difference of two singular terms.
Thus, for small values of $\lVert \bx - \bx' \rVert$, the kernel
$\psi_{\sigma}$ should be replaced by a suitable Taylor series
approximation.
Letting~$u = \lVert \bx - \bx' \rVert/\sqrt{2} \sigma(\bx)$, a
modest amount of algebra yields
\begin{equation}\label{eq:kerpsi_asymp}
\psi_\sigma(\bx,\bx') \approx 
     \hat{\vct{N}}(\vct{x}') \cdot \lp \vct{x}' -\vct{x} \rp  \lp
    \frac{ \frac{2}{3} - \frac{2u^2}{5} + \frac{u^4}{7} - \frac{u^6}{27} + \frac{u^8}{132}
+ \dots}{4\sqrt{2 \pi^3} \, \sigma(\bx)^3} \rp \, ,
\end{equation}
with an error of about $10^{-13}$ so long as $u <0.1$. For
$u \geq 0.1$, the loss of accuracy from catastrophic cancellation
in~\eqref{eq:kerpsi} is less than three digits of relative precision,
and therefore at least thirteen digits of relative accuracy are
obtained for any value of $u$ using standard double precision
arithmetic.

The gradient of $\Phi$ permits exactly the same decomposition.
Contributions from the far field can be obtained directly from the 
FMM and the near field can be computed directly.
The near field in \eqref{gradPhi} again involves the difference of
singular terms but can be replaced by the Taylor series approximation
\begin{multline}\label{eq:kerpsi_asymp2}
\frac{\partial \psi_\sigma(\bx,\bx')}{\partial x_i} \approx 
      -\hat{N}_i(\vct{x}') \lp
    \frac{\frac{2}{3} - \frac{2u^2}{5} + \frac{u^4}{7} -
      \frac{u^6}{27} + \frac{u^8}{132}+ \dots}{4\sqrt{2 \pi^3} \, \sigma(\bx)^3} \rp \\
+   \hat{\vct{N}}(\vct{x}') \cdot \lp \vct{x}' -\vct{x} \rp  \lp
\frac{e^{-u^2}}{2\sqrt{2 \pi^3} \, \sigma(\bx)^4} \rp
\frac{\partial \sigma(\bx)}{\partial x_i}  ,
\end{multline}
where~$\hat{N}_i(\vct{x}')$ denotes the $i$th component
of~$\hat{\vct{N}}(\vct{x}')$.

\subsection{Surface refinement}
\label{sec:refinement}

Since, as noted above, the level surface $\Gamma_{1/2}$ is defined
by~\eqref{discrete_int}, it is a simple matter to refine the
discretization so that the charts converge to $\Gamma_{1/2}$ with
$p$th-order accuracy.  One simply splits any skeleton
triangle~$\vct{T}^j$ where refinement is desired into four
subtriangles, adding the midpoints of each side as new vertices. In
general, this procedure will break the earlier assumption of
conformity of the skeleton mesh, but it is easy to check that the same
refinement of the pseudonormal vector field remains continuous.  One
can then construct the pseudonormal vector fields on the subtriangles
as above and solve the nonlinear equation~\eqref{eq:newton} for each
of the new quadrature nodes. Determining which sections of the
geometry~$\Gamma^j$ need to be refined is up to the user, and standard
\emph{a posteriori} error estimates
from adaptive interpolation or adaptive integration work well.

\begin{remark}
  The method described in this paper is based on the fundamental
  premise that the skeleton mesh is well enough resolved that the
  ``nearby" $C^\infty$ surface we construct can serve the purposes of
  the subsequent simulation.  If $\Gamma_{1/2}$ is
  \emph{unsatisfactory}, additional tools would be required that
  enable modification of the input triangulation.  This is beyond the
  scope of the present work.
\end{remark}

%
%
%

\section{Numerical examples}
\label{sec:examples}

In this section we provide several numerical examples demonstrating
the behavior and computational efficiency of our algorithm in
converting flat skeleton triangulations into higher order surfaces. In
addition to computational scaling results, we also present results of
computing a Gauss flux integral to estimate how \emph{water-tight} the
smooth surface is. To this end, let~$\bx_0$ denote some point in the
interior of a region~$\Omega$ with boundary~$\Gamma$. It is easy to see
that by straightforward application of the divergence theorem, and
using the fact that the Green's function of the Laplace operator
is~$1/4\pi r$:
\begin{equation}
  \begin{aligned}
    1 &= \int_\Omega \delta(\bx_0 - \bx) \, dv(\bx) \\
    &= -\int_\Omega \Delta \frac{1}{4\pi \lVert \bx_0 - \bx \rVert} \,
    dv(\bx) \\
    &= -\int_\Omega \nabla \cdot \nabla \frac{1}{4\pi \lVert \bx_0 - \bx \rVert} 
    \, dv(\bx) \\
    &= -\int_\Gamma \hat{\vct{N}}(\bx) \cdot \nabla \frac{1}{4\pi \lVert
      \bx_0
      - \bx \rVert} 
    \, dv(\bx).
  \end{aligned}
\end{equation}
The last expression is a measure of the flux through the
surface~$\Gamma$, and we therefore estimate the deviation of the
output of our algorithm from a closed surface as~$\epsilon_0$:
\begin{equation}
  \epsilon_0 = 1 + \int_{\Gamma_{1/2}} \hat{\vct{N}}(\bx) \cdot \nabla
  \frac{1}{4\pi \lVert  \bx_0 - \bx \rVert}   \, dv(\bx).
\end{equation}

In addition to reporting the value of~$\epsilon_0$, tables containing
results for the following numerical experiments also contain the data:
\begin{itemize}
\item $M$, the number of skeleton triangles,
\item $q$, the order of Vioreanu-Rokhlin quadrature along~$S$,
\item $n_q$, the number of quadrature nodes on each skeleton
  triangle,
\item $p$, the order of the discretization of~$\Gamma_{1/2}$,
\item $n_p$, the number of nodes on each panel~$\Gamma^j$,
\item $k$, number of Newton iterations,
\item $\epsilon_N$, max point-wise error in $h$
  after Newton iterations, as measured by
  \begin{equation}
    \epsilon_N = \max_{j,i} \left\vert h^{(k)}_{ji} - h^{(k-1)}_{ji} \right\vert,
  \end{equation}
\item $\epsilon_0$, error in flux integral, described above,
\item $\tfmm$, the time required for a single FMM call, i.e. time
  for one Newton iteration,
\end{itemize}
In each problem, the Newton iterations were run to a tolerance
of~$10^{-12}$.
In practice, the convergence of each Newton iteration for each
discretization point along~$\Gamma_{1/2}$ is independent, and
therefore points (targets) that have already converged can be removed
from the target list for subsequent FMM calls. A modest reduction in
computational cost can be obtained by this optimization, but the
overall scheme remains linear in cost and dominated by the initial
Newton iterations.

Lastly, the algorithm was implemented 
in Fortran 77/95, and compiled with the Intel Fortran Compiler
2019. Examples were run on a
workstation with 32 Intel Xeon Gold 6130 cores at 2.1GHz with 512GB of
shared memory;
modest multicore acceleration of the FMM was done using OpenMP
directives.
Plots of 3D images were created in Paraview~\cite{ahrens2005paraview}.

\subsection{Basic surface construction}

In this example, we merely show the results of our algorithm when
applied to a flat triangulation of a smooth surface, namely that of a
cube with rounded edges and corners. The geometry was constructed in
FreeCAD~\cite{freecad}, exported as a \texttt{.step} file, and then
imported and meshed in GiD v13.0.4 on Linux~\cite{gid}.
Two skeleton meshes are used in the example: one with 5736 flat
triangles and a refined version with 21,852 flat triangles.
The skeleton meshes and
smoothed surfaces are shown in Figure~\ref{fig:rcube}.
Data for this example is contained in Table~\ref{tab:rcube}.
The original un-rounded cube has one vertex at $(0,0,0)$ and the
opposite vertex at~$(10,10,10)$. In order to test the Gauss flux
integral, a point-source was placed at $(4.5,4.5,5)$.

Table~\ref{tab:rcube_conv} contains results from a refinement experiment.
In this study, each of the skeleton triangles is refined into four 
sub-triangles, and the smoothing procedure is carried out again. While this
refinement does not alter the skeleton surface, it does scale the width
of the smoothing kernel~$\sigma$, as well as result in a more resolved
approximation of the smoothed surface (as shown by the convergence 
of~$\epsilon_0$). Near machine precision water-tightness is relatively easy
to achieve.

\afterpage{

\begin{table}[t]
  \centering
  \caption{Results for smoothing a triangulation of a rounded cube.}
  \label{tab:rcube}
    \begin{tabular}{|c|c|c|c|c|c|c|c|c|c|c|c|}
      \hline
      $M$ &$q$ & $n_q$ & $Mn_q$ & $p$ &$n_p$ & $Mn_p$ & $k$ 
      &$\epsilon_N$ &$\epsilon_0$  & $\tfmm$ \\ \hline
      5736 & 4  & 15 & 86,040  & 4 & 15 & 86,040 & 4 & 1.0 E-12 
          & 1.1 E-07 & 5.9 E+00 \\
      5736 & 8  & 45 & 258,120 & 4 & 15 & 86,040 & 4 & 1.0 E-12 
          &  1.1 E-07 & 1.1 E+01 \\
      5736 & 12 & 91 & 521,976 & 4 & 15 & 86,040 & 4 & 1.0 E-12
          & 1.1 E-07 & 1.9 E+01 \\
      \hline
       5736 & 4 & 15 & 86,040  & 8 & 45 & 258,120 & 4 & 1.0 E-12 
          & 4.7 E-13 & 1.5 E+01 \\
       5736 & 8 & 45 & 258,120 & 8 & 45 & 258,120 & 4 & 1.1 E-12 
          & 4.8 E-13 & 2.9 E+01\\
      \hline 
    \end{tabular}
\end{table}

\begin{table}[t]
  \centering
  \caption{Convergence under refinement for smoothing 
  a triangulation of a rounded cube.}
  \label{tab:rcube_conv}
    \begin{tabular}{|c|c|c|c|c|c|c|c|c|c|c|c|}
      \hline
      $M$ &$q$ & $n_q$ & $Mn_q$ & $p$ &$n_p$ & $Mn_p$ & $k$ 
      &$\epsilon_N$ &$\epsilon_0$  & $\tfmm$ \\ \hline
      5736 & 4  & 15 & 86,040  & 4 & 15 & 86,040 & 4 & 1.0 E-12 
          & 1.1 E-07 & 5.9 E+00 \\
      22,944 & 4  & 15 & 344,160 & 4 & 15 & 344,160 & 3 & 1.0 E-12 
          & 1.9 E-10 & 1.9 E+01 \\
      91,776 & 4  & 15 & 1,376,640 & 4 & 15 & 1,376,640 & 2 & 1.0 E-12 
          & 4.6 E-13 & 7.2 E+01 \\
      \hline 
    \end{tabular}
\end{table}

\begin{figure}[b]
  \centering
  \begin{subfigure}[b]{.4\linewidth}
    \centering
    \includegraphics[width=.95\linewidth]{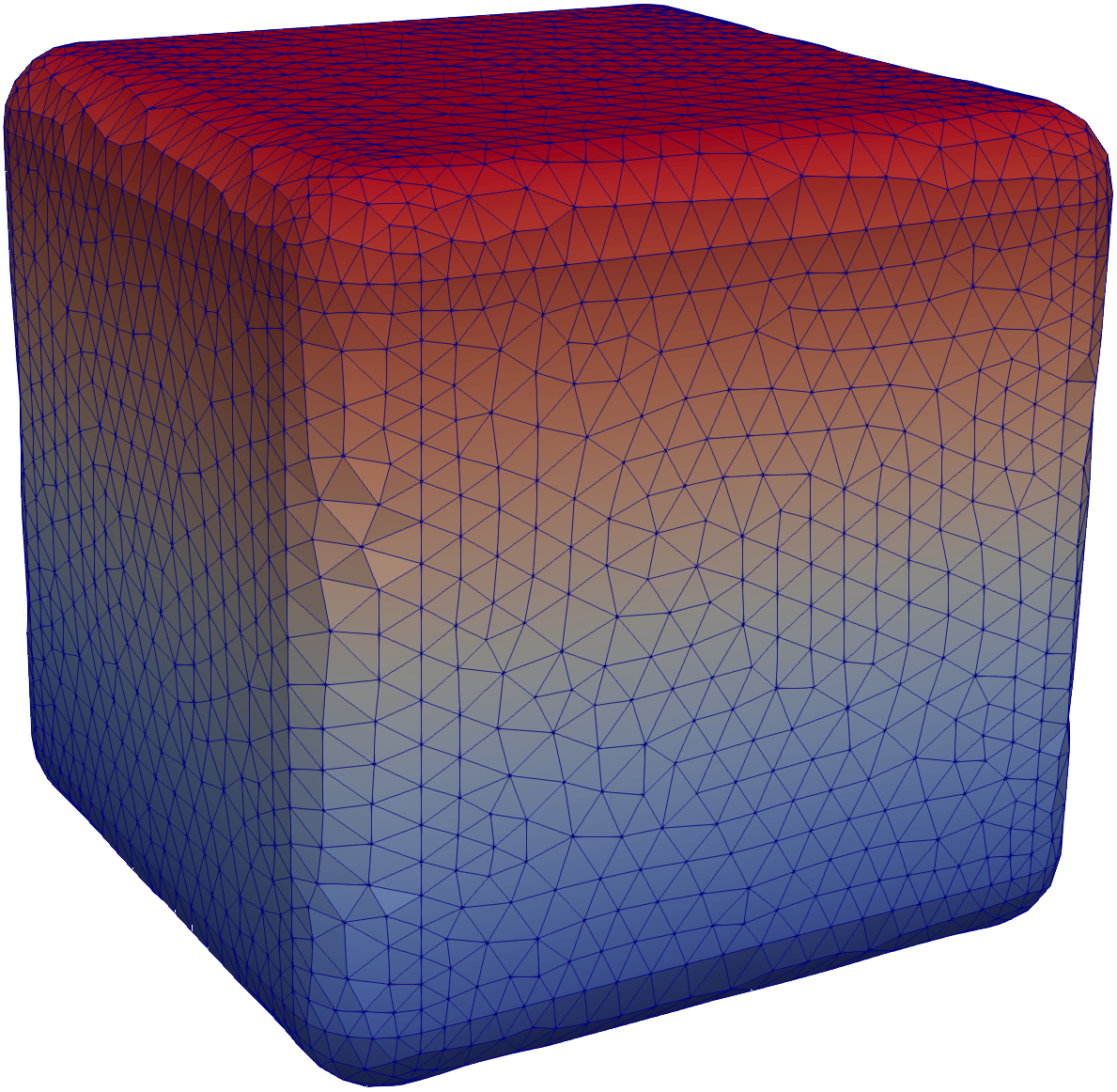}
    \caption{Skeleton mesh.}
    \label{fig:rcube5736-skel}
  \end{subfigure}
  \hfill
  \begin{subfigure}[b]{.4\linewidth}
    \centering
    \includegraphics[width=.95\linewidth]{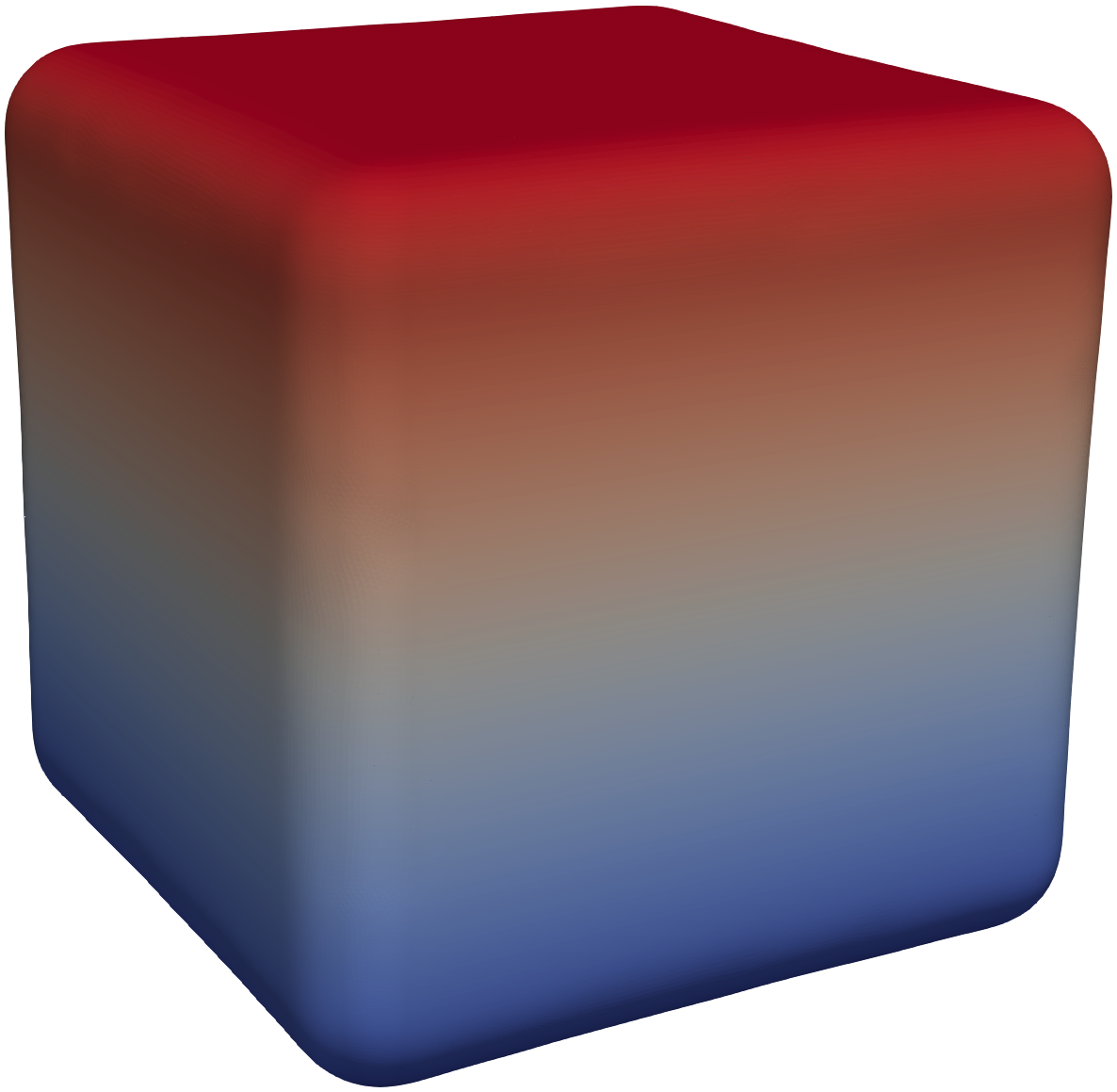}
    \caption{Smoothed surface from~\ref{fig:rcube5736-skel}.}
    \label{fig:rcube5736-smoo}
  \end{subfigure} \\
  \vspace{\baselineskip}
\begin{subfigure}[b]{.4\linewidth}
    \centering
    \includegraphics[width=.95\linewidth]{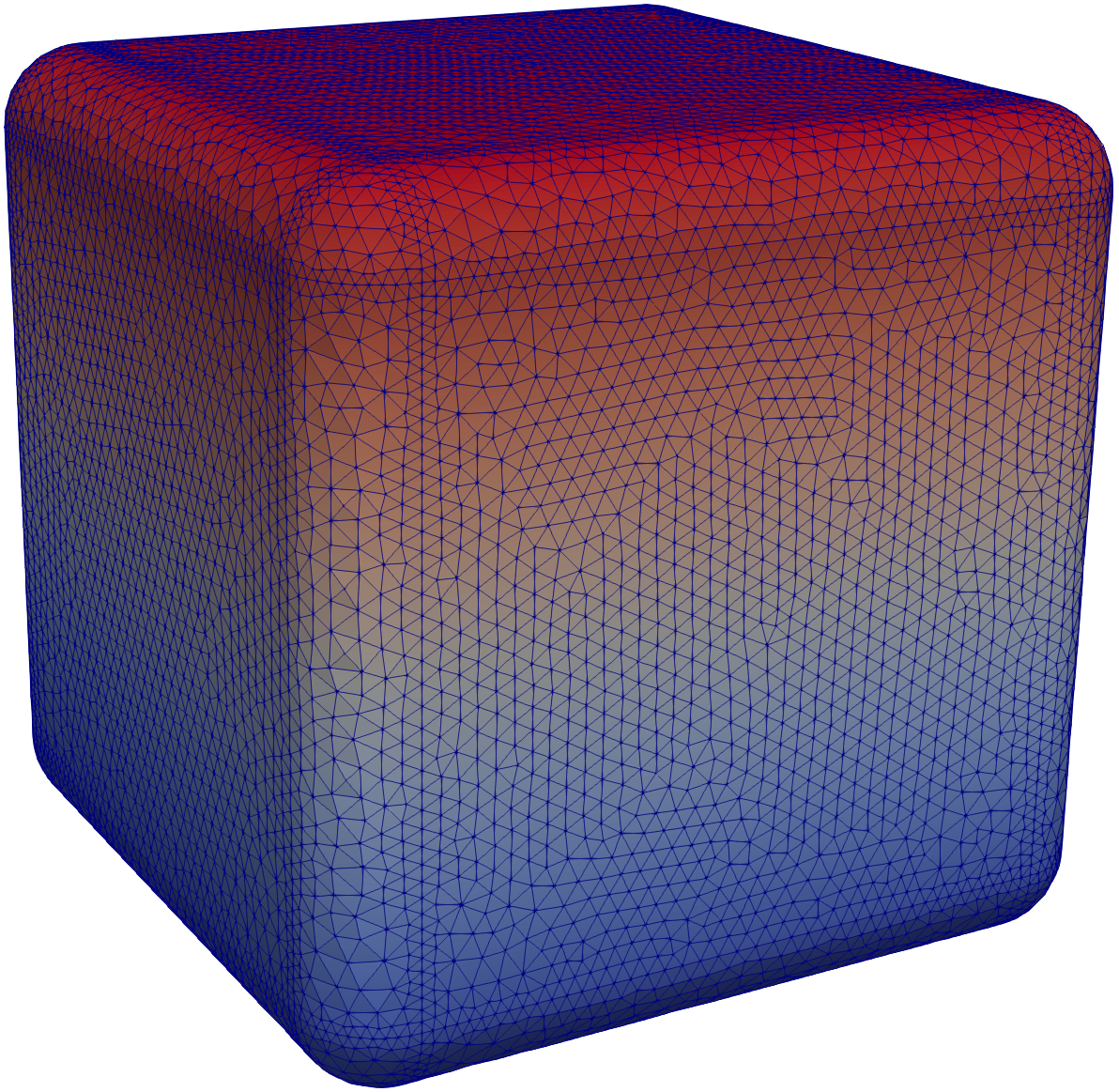}
    \caption{Skeleton mesh.}
    \label{fig:rcube21852-skel}
  \end{subfigure}
  \hfill
  \begin{subfigure}[b]{.4\linewidth}
    \centering
    \includegraphics[width=.95\linewidth]{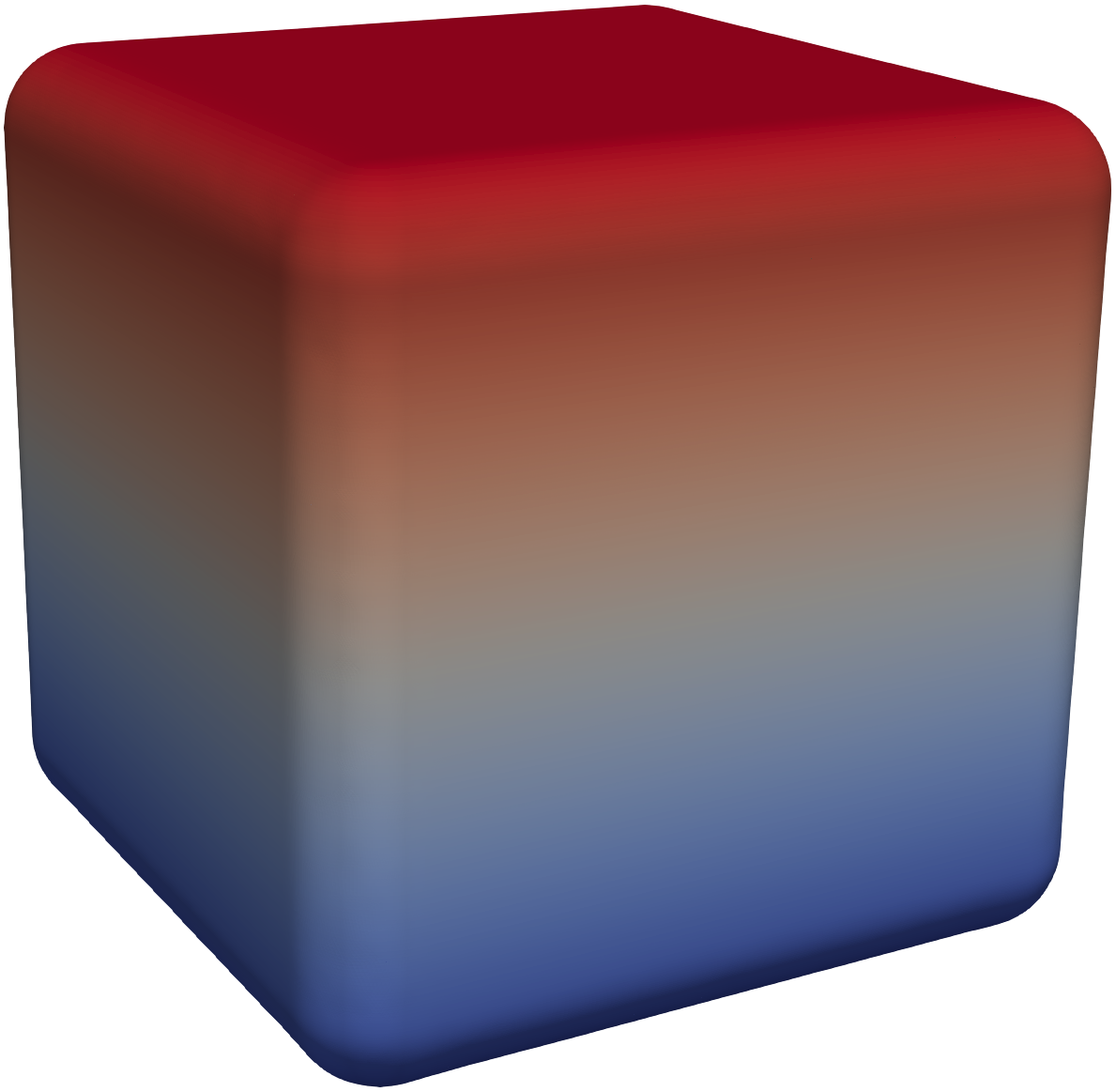}
    \caption{Smoothed surface from~\ref{fig:rcube21852-skel}.}
    \label{fig:rcube21852-smoo}
  \end{subfigure}  
  \caption{Skeleton and smooth surfaces for a rounded cube. The
    skeleton mesh in~\ref{fig:rcube5736-skel} contains 5736 flat
    triangles, and the resulting smoothed surface
    in~\ref{fig:rcube5736-smoo} is piecewise 4th-order.  The skeleton
    mesh in~\ref{fig:rcube21852-skel} contains 21,852 flat triangles,
    and the resulting smoothed surface in~\ref{fig:rcube21852-skel} is
    piecewise 4th-order. In each case, the integral defining the level
    set was discretized using 4th-order Vioreanu-Rokhlin quadratures.}
  \label{fig:rcube}
\end{figure}
\clearpage}

\subsection{Recovery of a sphere}
\label{sec:sphere}

\begin{figure}[t]
  \centering
  \begin{subfigure}[b]{.4\linewidth}
    \centering
    \includegraphics[width=.95\linewidth]{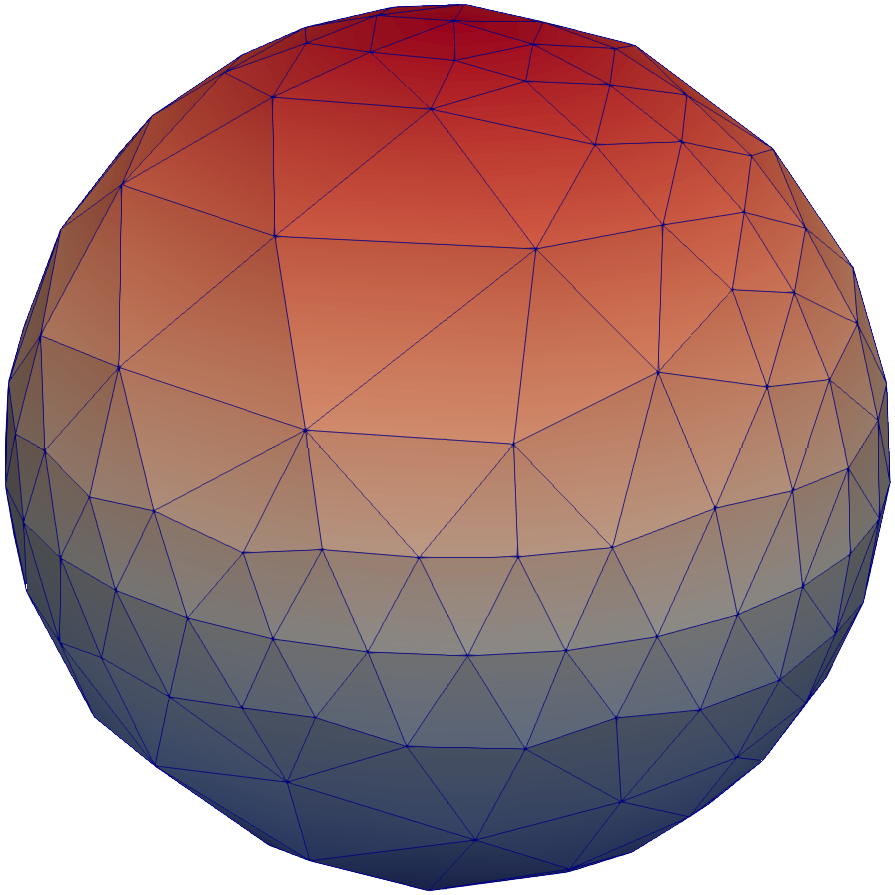}
    \caption{The skeleton mesh.}
  \end{subfigure}
  \hfill
  \begin{subfigure}[b]{.4\linewidth}
    \centering
    \includegraphics[width=.95\linewidth]{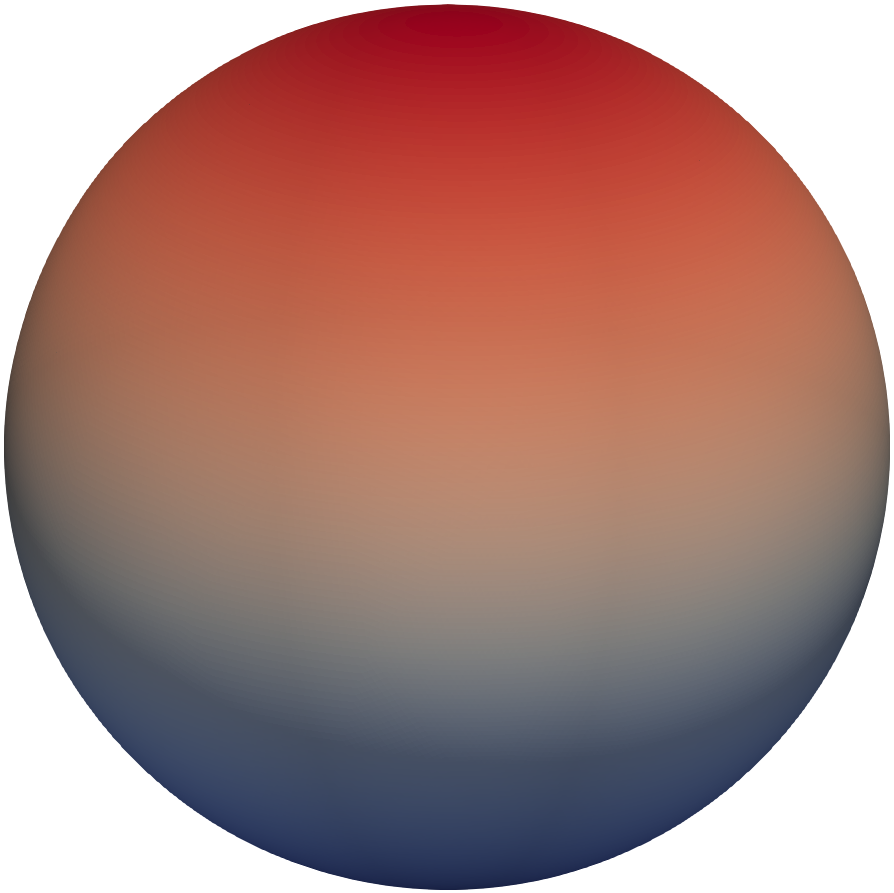}
    \caption{The smoothed surface.}
  \end{subfigure}
  \caption{Skeleton and smooth surfaces for a sphere. The skeleton
    mesh contains 416 triangles. The skeleton surface was discretized
    using 4th-order quadratures, and the smoothed surface is piecewise
    4th-order smooth.}
  \label{fig:sphere}
\end{figure}

In this example, we obtain a smooth surface from a rough flat
triangulations of a sphere. Figure~\ref{fig:sphere} illustrates an
example skeleton mesh and the smooth surface obtained using the
algorithm of this work. As an additional measure, we also report the
average of the norm of the discretization points (i.e. distance from
the origin) on the smooth surface, as well as the standard deviation:
\begin{equation}
\langle \vct{x} \rangle = \frac{1}{M n_p} \sum_{i=1}^M \sum_{j=1}^{n_p}
\left\Vert \bx_{ij} \right\Vert, \qquad \textnormal{std}(\vct{x}) =
\sqrt{\frac{1}{Mn_p} 
  \sum_{i=1}^M \sum_{j=1}^{n_p} \lp \left\Vert \vct{x}_{ij} \right\Vert
  - \langle \vct{x}
  \rangle \rp^2
}
\end{equation}
Table~\ref{tab:sphere} contains results. In order to check the Gauss
integral, a point source was placed at~$(0.1, 0, 0)$. The triangulation was
obtained directly from GiD~\cite{gid}, and little effort was made to
ensure that it was regular. Clearly there is a modest amount of
adaptive refinement near the edges of the CAD surfaces, but such
features are not visible in the smoothed geometry. Lastly, note that
Table~\ref{tab:sphere} contains the dual information as
Table~\ref{tab:rcube}: $p$ is varied for fixed~$q$, instead of
vice versa. In summary, the smooth surface is a smooth sphere-like
object that deviates from true sphere by roughly 0.01, albeit
with a radius not equal to one. This is quite good, as the original
skeleton mesh was not particularly fine.

\begin{table}[b]
  \centering
  \caption{Results for smoothing a triangulation of a sphere.}
  \label{tab:sphere}
  \resizebox{\textwidth}{!} {
    \begin{tabular}{|c|c|c|c|c|c|c|c|c|c|c|c|c|}
      \hline
      $M$ &$q$ & $n_q$ & $Mn_q$ & $p$ &$n_p$ & $Mn_p$
      & $k$ & $\epsilon_N$
      &$\epsilon_0$  & $\tfmm$ & $\langle \vct{x} \rangle$
      & std$(\vct{x})$\\ \hline 
      416 & 4 & 15 & 6240 & 4 & 15 & 6240 & 4 & 4.2 E-13
          & 9.6 E-08 & 1.1 E+00 & 0.96473 & 1.4 E-02 \\
      416 & 4 & 15 & 6240 & 8 & 45 & 18,720 & 4 & 9.6 E-13
         & 6.2 E-13 & 3.1 E+00 & 0.96375 & 1.4 E-02 \\
      416 & 4 & 15 & 6240 & 12 & 91 & 37,856 & 4 & 1.0 E-12
         & 3.0 E-15 & 6.4 E+00 & 0.96375 & 1.4 E-02 \\
      \hline
      416 & 8 & 45 & 18,720 & 4 & 15 & 6240 & 4 & 9.8 E-13
          & 9.6 E-07 & 3.5 E+00 & 0.96373 & 1.4 E-02 \\
      416 & 8 & 45 & 18,720 & 8 & 45 & 18,720 & 4 & 6.9 E-13
         & 3.9 E-13 & 9.1 E+00 & 0.96375 & 1.4 E-02 \\
      416 & 8 & 45 & 18,720 & 12 & 91 & 37,856 & 4 & 1.0 E-12
         & 2.2 E-15 & 1.9 E+01 & 0.96375 & 1.4 E-02 \\
      \hline 
    \end{tabular}
    }
\end{table}

\subsection{Quadratic skeletons}
\label{sec:quadratic}

A torus is one of the smooth primitives in the software
GiD~\cite{gid}, and therefore a straightforward test.
In this example, we show the remarkable increase in quality of the
smooth surface when the skeleton mesh consists of quadratic patches
(i.e. 2nd-order curvilinear triangles) instead of flat triangles. In
order to clearly demonstrate the benefit, we have set~$\lambda=10$
in~\eqref{sigma2}
(the scaling parameter for the width of the convolution kernel).
This narrows the width of the convolving kernel, which allows for
preservation of more fine-scale features of the actual skeleton mesh.
We only present images as qualitative results in this
case, as the various convergence results are similar to the previous
two examples. It is worth pointing out that quadratic triangles can
easily be shown to form an exact water-tight surface, and therefore
the construction of our pseudonormal vector field is analogous to the
case in which the skeleton mesh consists of flat triangles.

\subsection{Large-scale structures}
\label{sec:large}

\afterpage{
\begin{figure}[t]
  \centering
  \begin{subfigure}[b]{.45\linewidth}
    \centering
    \includegraphics[width=.95\linewidth]{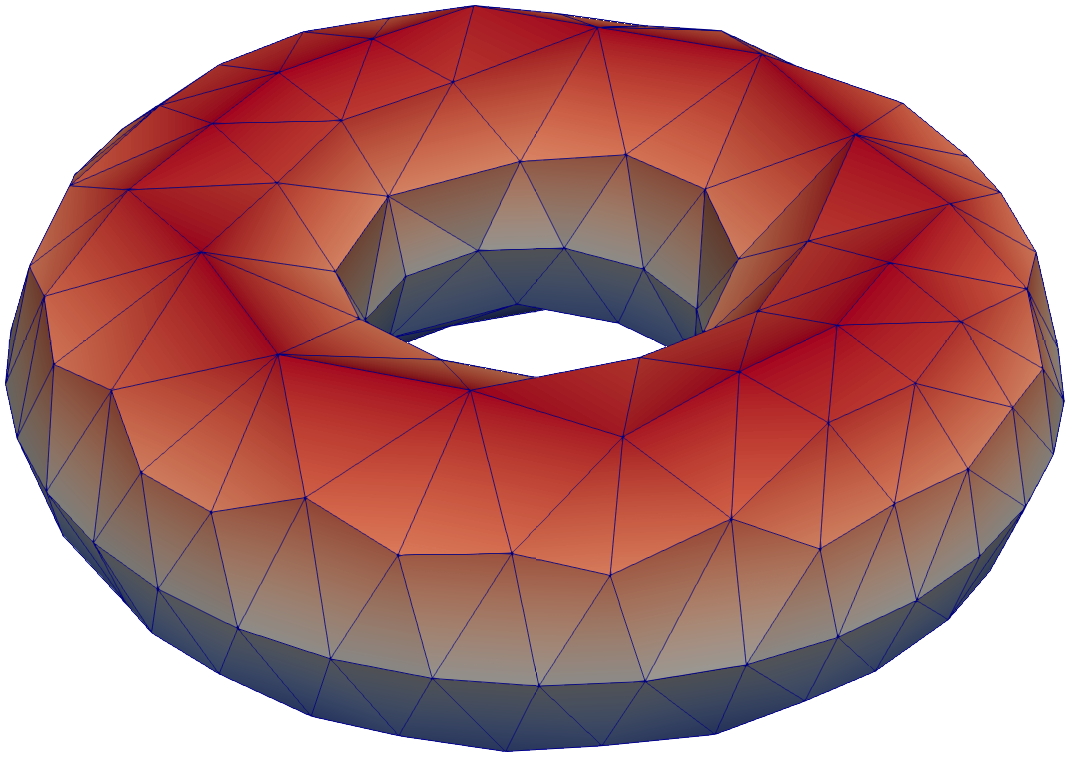}
    \caption{A flat skeleton.}
  \end{subfigure}
  \hfill
  \begin{subfigure}[b]{.45\linewidth}
    \centering
    \includegraphics[width=.95\linewidth]{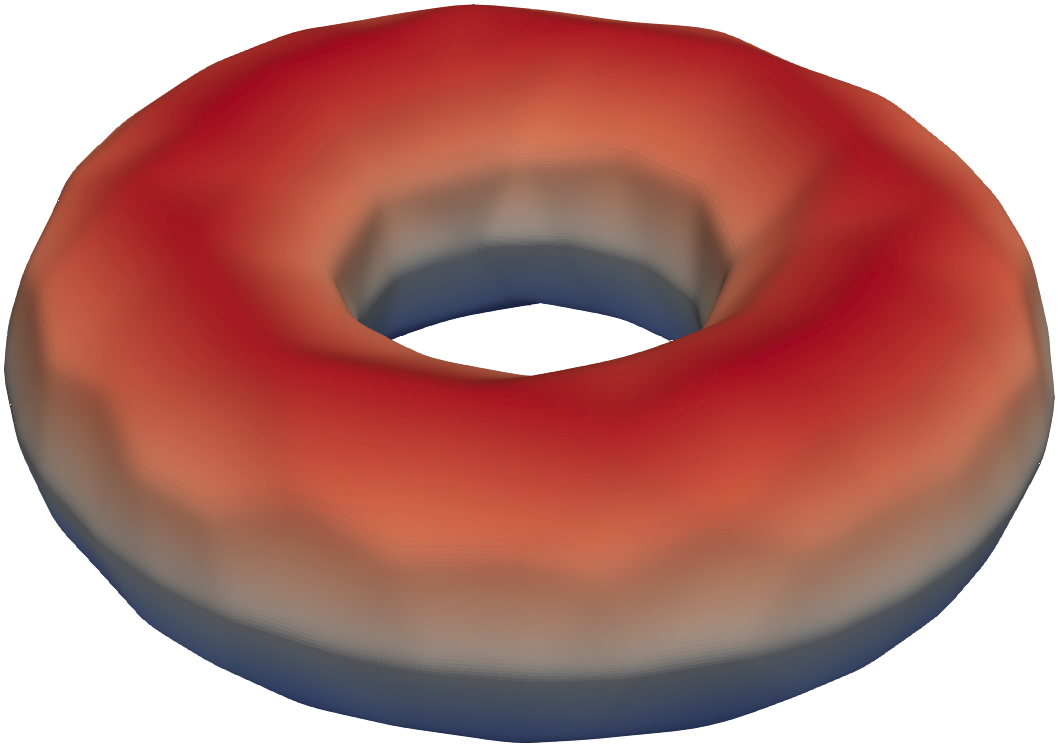}
    \caption{$\Phi_{1/2}$ from flat skeleton.}
  \end{subfigure}\\
  \vspace{\baselineskip}
  \begin{subfigure}[b]{.45\linewidth}
    \centering
    \includegraphics[width=.95\linewidth]{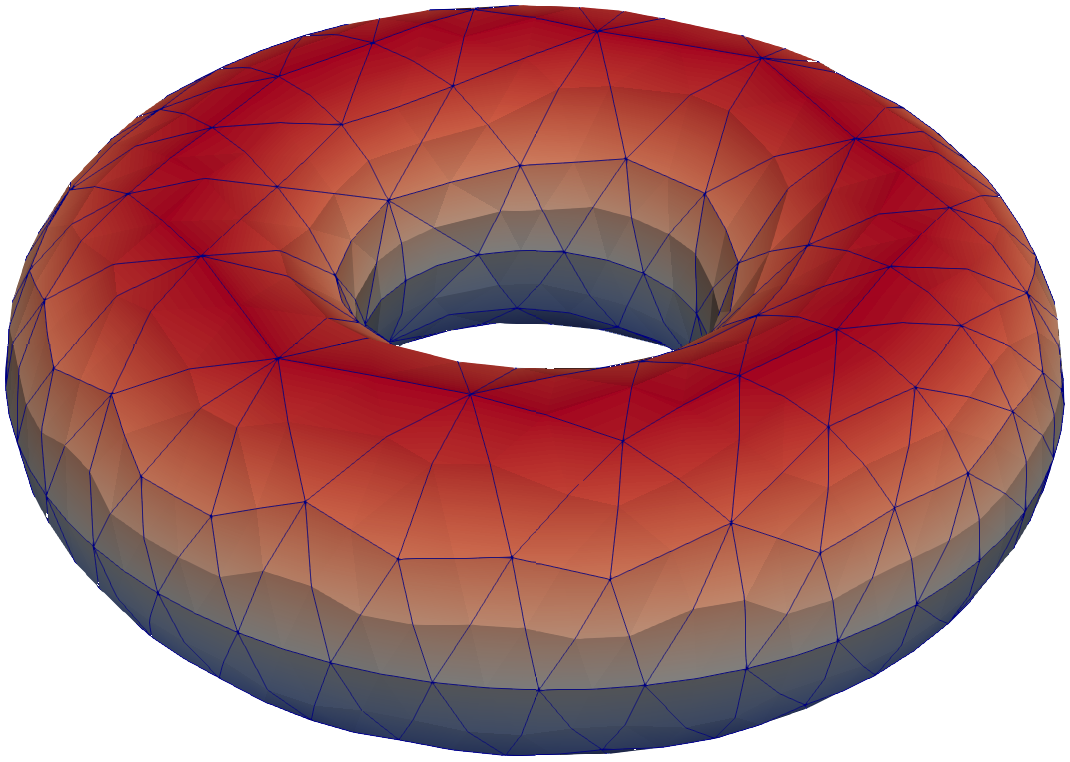}
    \caption{A quadratic skeleton.}
  \end{subfigure}
  \hfill
  \begin{subfigure}[b]{.45\linewidth}
    \centering
    \includegraphics[width=.95\linewidth]{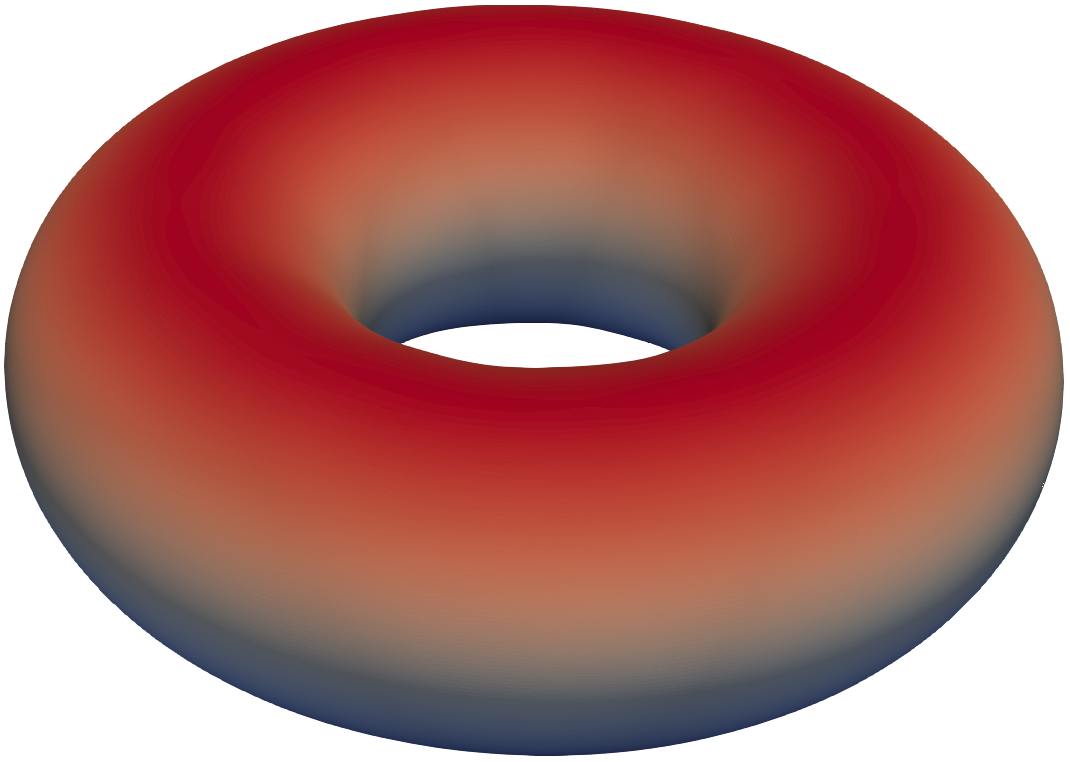}
    \caption{$\Phi_{1/2}$ from quadratic skeleton.}
  \end{subfigure}
  \caption{A qualitative comparison of the improvement when using
    quadratic skeleton meshes instead of flat ones. In each case, the
    skeleton surface was discretized using 12th-order quadratures, and
    the smoothed surface is piecewise smooth to 12th-order. In both
    examples, the Gauss integral was accurate to approximately 9 digits.}
\end{figure}
\clearpage}

In this section we demonstrate the performance of our algorithm on a
relatively large-scale structure with highly multiscale
features. Figure~\ref{fig:large} contains images of the skeleton mesh
and resulting smoothed surface for a mock-up of an A380 passenger
aircraft with small antennas mounted on the top. The geometry was
designed and meshed in GiD v13.0.4. The skeleton mesh
consists of 31,336 quadratic curvilinear triangles. Each skeleton
patch was discretized using a 4th order Vioreanu-Rokhlin quadrature,
we set~$\lambda=10$, 
and the smooth surface was computed to 4th order as well.
Despite the existence of edges on the multiscale antennas, as seen in
the figures, the algorithm does a qualitatively good job of obtaining
a nearby smooth surface without any noticeable ringing
artifacts. Convergence in the divergence test is of high-order, and
the FMM scales linearly (as expected). This data is presented in
Table~\ref{tab:large_conv}.

\afterpage{
\begin{figure}[t]
  \centering
  \begin{subfigure}[b]{.55\linewidth}
    \centering
    \includegraphics[width=.95\linewidth]{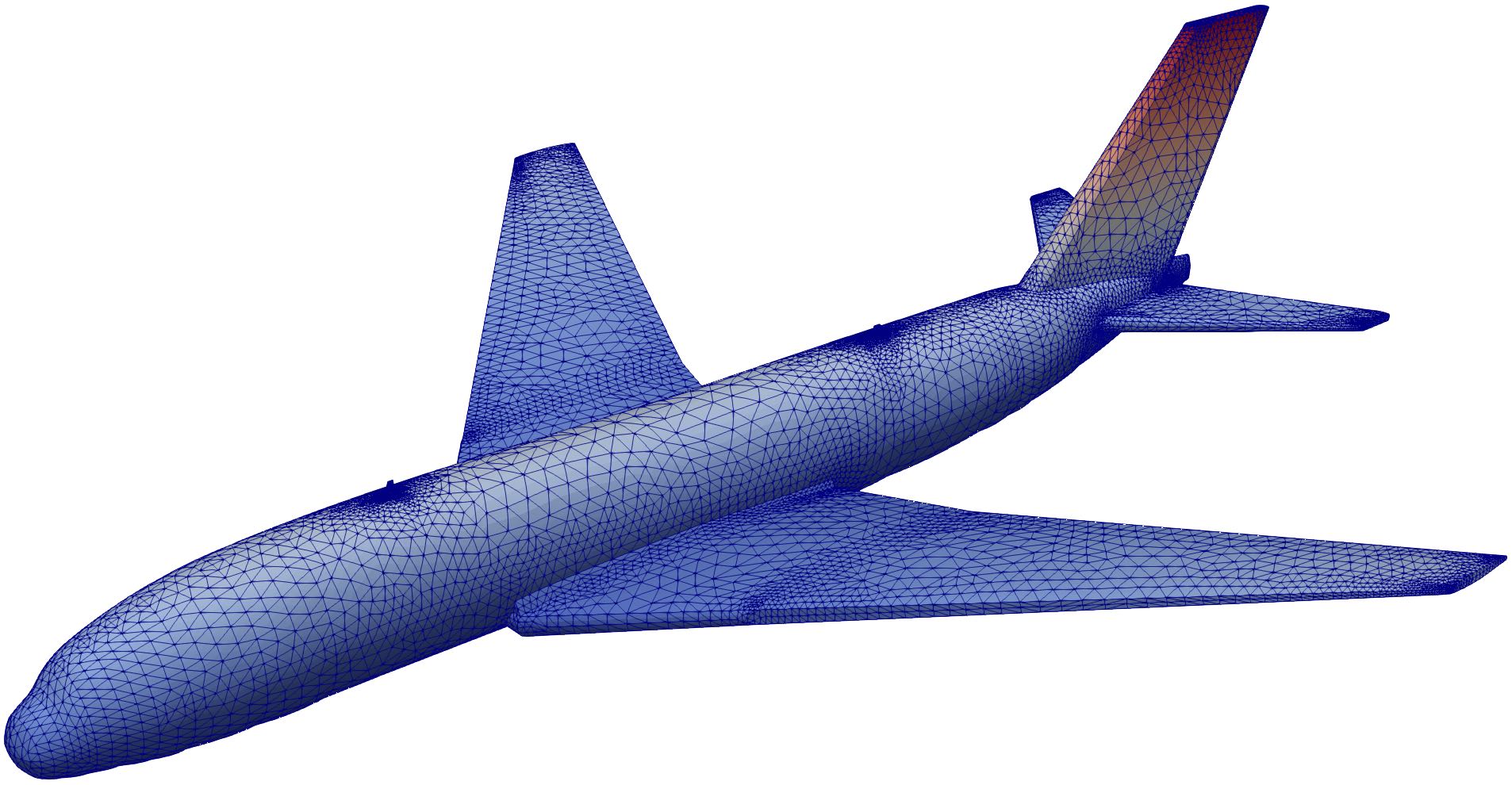}
    \caption{The skeleton mesh.}
  \end{subfigure}
  \hfill
  \begin{subfigure}[b]{.4\linewidth}
    \centering
    \includegraphics[width=.95\linewidth]{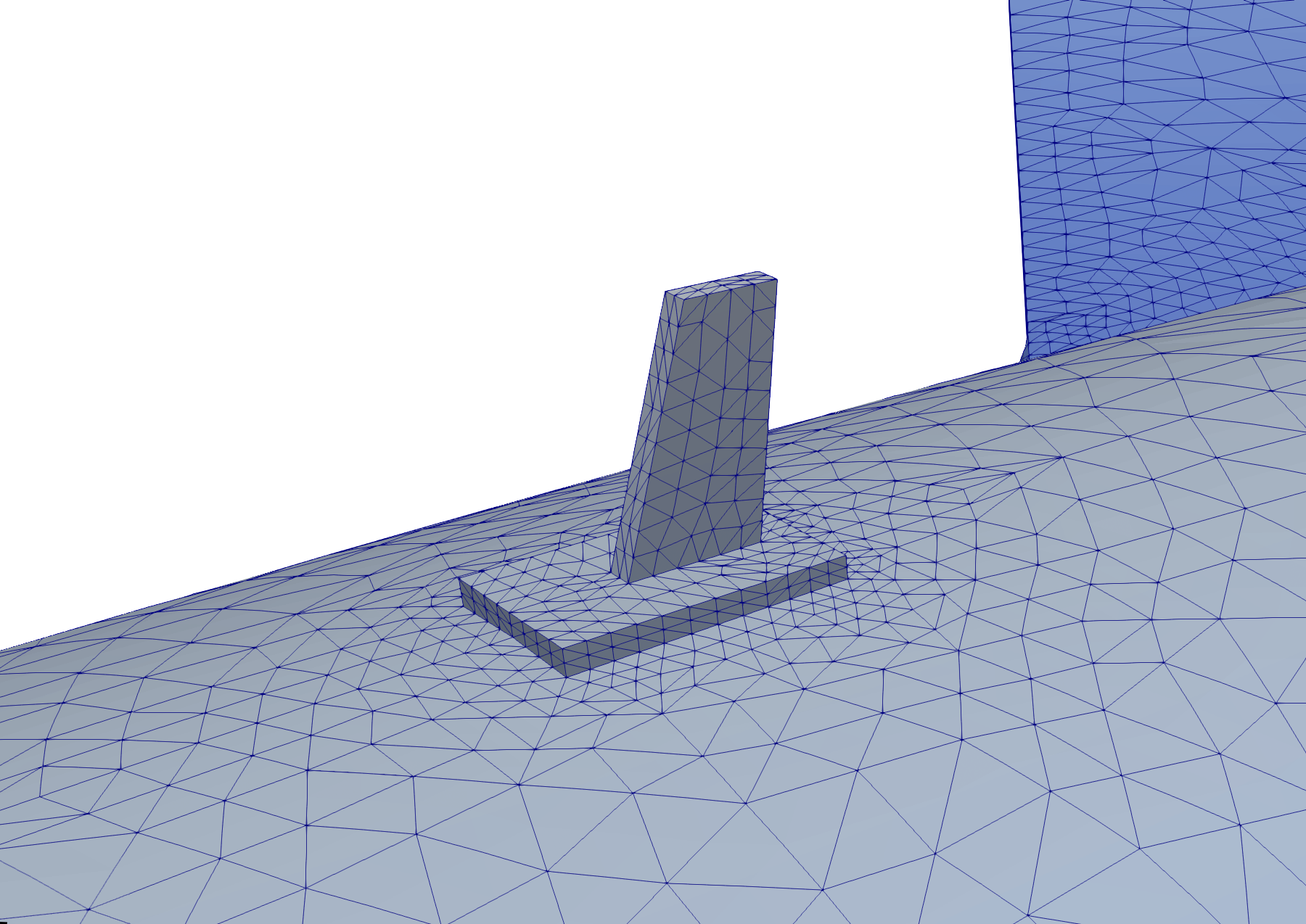}
    \caption{Multiscale antenna.}
  \end{subfigure}\\
  \vspace{\baselineskip}
  \begin{subfigure}[b]{.55\linewidth}
    \centering
    \includegraphics[width=.95\linewidth]{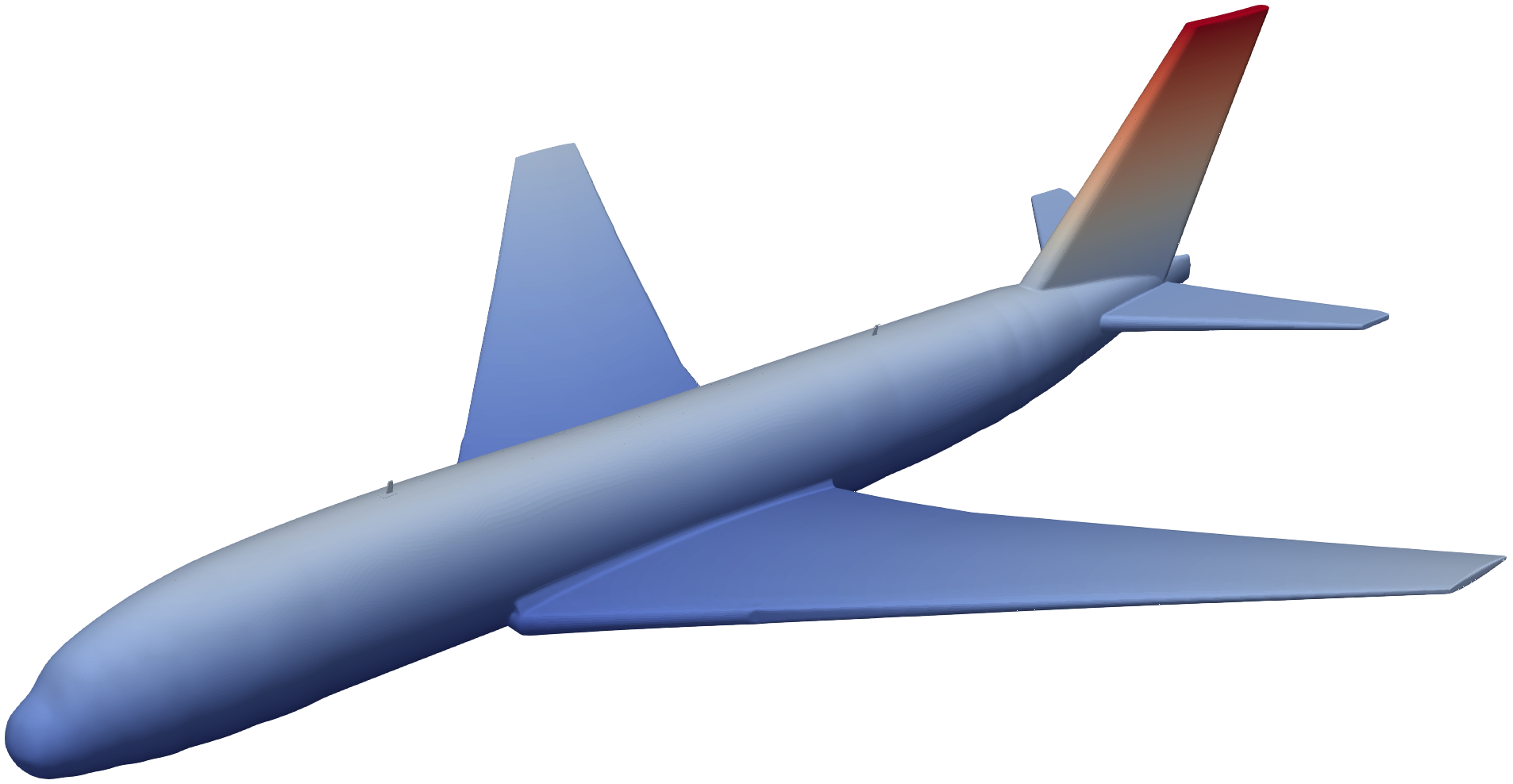}
    \caption{The smoothed surface.}
  \end{subfigure}
  \hfill
  \begin{subfigure}[b]{.4\linewidth}
    \centering
    \includegraphics[width=.95\linewidth]{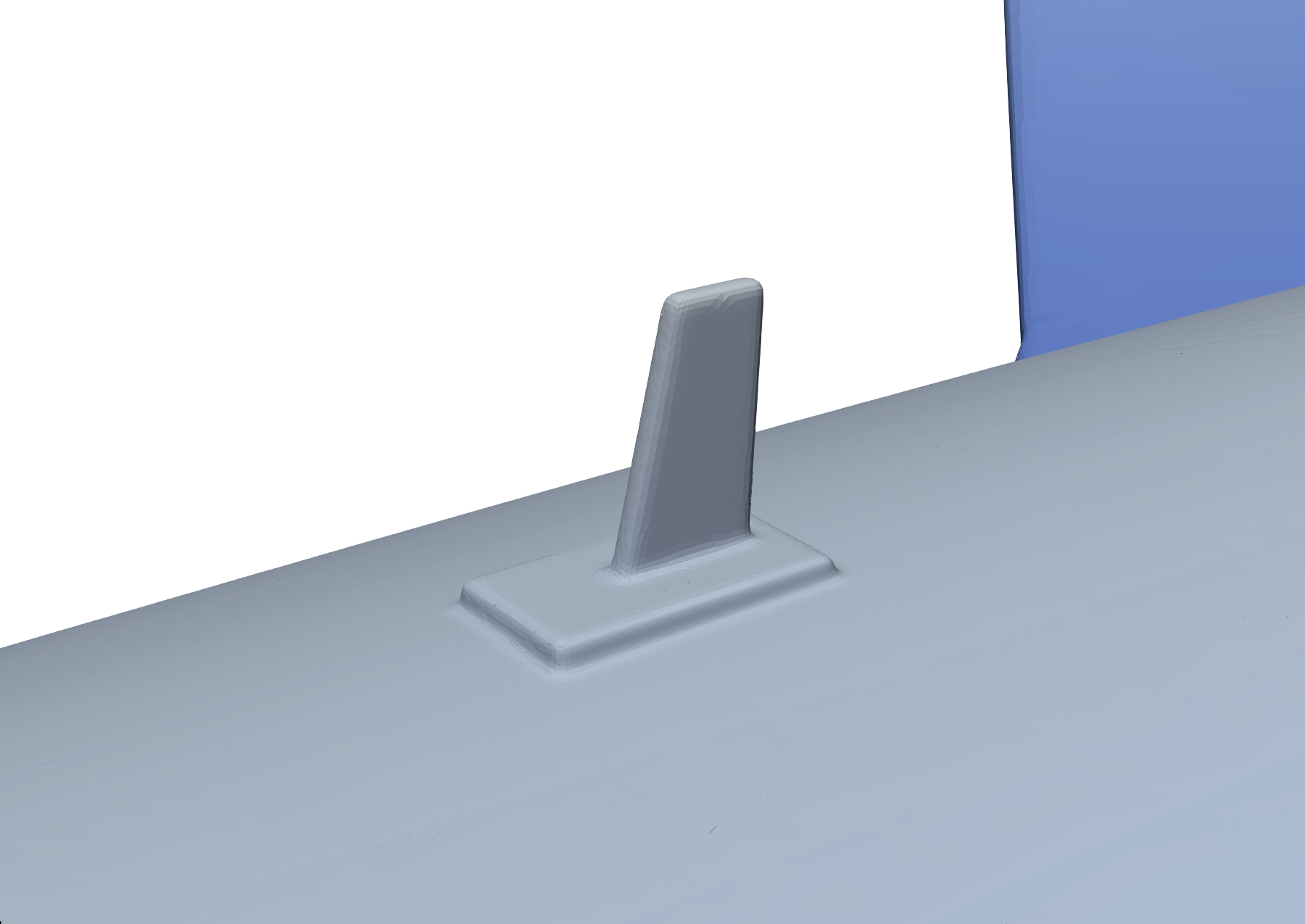}
    \caption{The smoothed antenna.}
  \end{subfigure}
  \caption{Quadratic skeleton mesh and smooth surfaces for a model
    passenger airplane. The skeleton mesh contains 31,336 quadratic
    triangles, and the smoothed surface is piecewise 4th-order
    smooth.}
  \label{fig:large}
\end{figure}

\begin{table}[b]
  \centering
  \caption{Convergence under refinement for smoothing 
  a quadratic triangulation of a passenger airplane.}
  \label{tab:large_conv}
    \begin{tabular}{|c|c|c|c|c|c|c|c|c|c|c|c|}
      \hline
      $M$ &$q$ & $n_q$ & $Mn_q$ & $p$ &$n_p$ & $Mn_p$ & $k$ 
      &$\epsilon_N$ &$\epsilon_0$  & $\tfmm$ \\ \hline
      31,336 & 4  & 15  & 470,040  & 4 & 15 & 470,040 & 4 & 1.0 E-12 
          & 1.5 E-06 & 1.8 E+01 \\
      125,344 & 4  & 15  & 1,880,160  & 4 & 15 & 1,880,160 & 4 & 1.0 E-12 
          & 2.1 E-08 & 6.3 E+01 \\
      501,376 & 4  & 15  & 7,520,640  & 4 & 15 & 7,520,640 & 4 & 1.0 E-12 
          & 3.6 E-10 & 2.5 E+02 \\
      \hline 
    \end{tabular}
\end{table}

\clearpage}

\subsection{Modes of failure}
\label{sec:failure}

There is one main regime in which our algorithm does not produce
qualitatively acceptable results: when the skeleton mesh is
too coarse. On the one hand, since our algorithm is
attempting to recover a smooth surface nearby to the skeleton mesh,
if the skeleton mesh is very coarse the argument could be made that it
is not a good approximation of an underlying smooth
surface. On the other hand, there are many geometries in which a
coarse mesh is an exact description of the geometry, for example 12
flat triangles exactly describe the surface of a cube.

\begin{figure}[t]
  \centering
  \begin{subfigure}[b]{.45\linewidth}
    \centering
    \includegraphics[width=.95\linewidth]{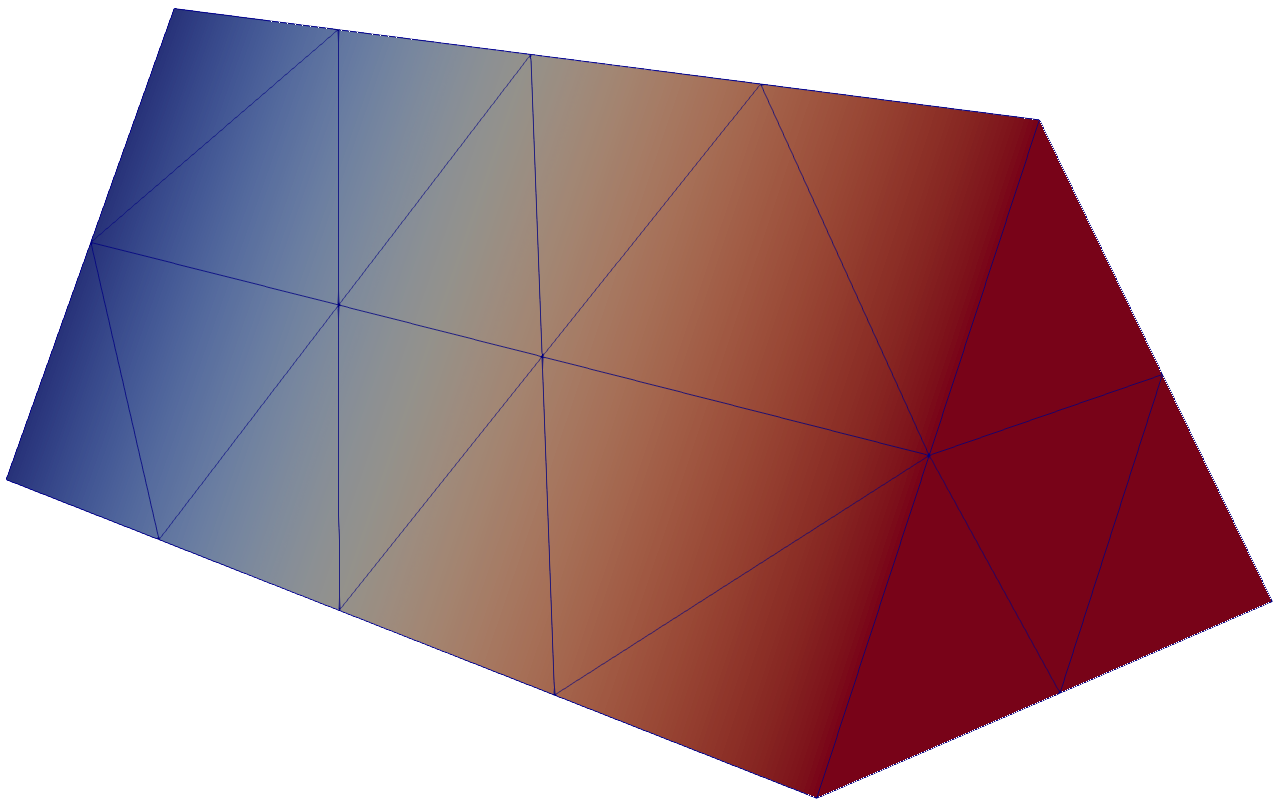}
    \caption{A very coarse skeleton mesh with 50 elements.}
  \end{subfigure}
  \hfill
  \begin{subfigure}[b]{.45\linewidth}
    \centering
    \includegraphics[width=.95\linewidth]{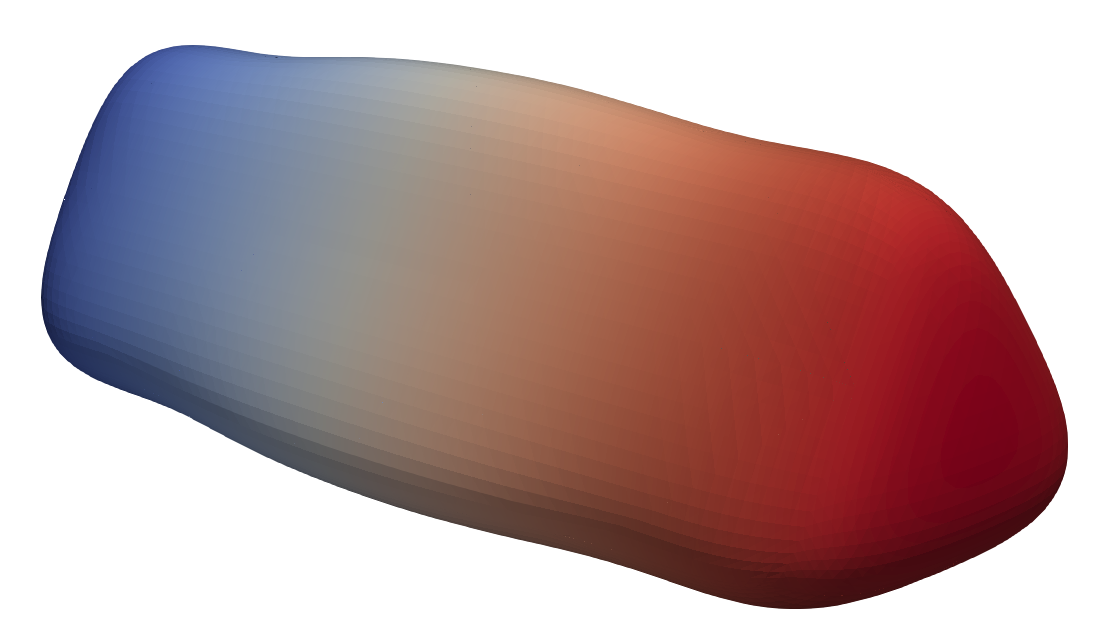}
    \caption{A naive smoothing of coarse mesh with~$\lambda =
      5$.}
    \label{fig:prism_ringing}
  \end{subfigure}\\
  \vspace{\baselineskip}
  \begin{subfigure}[b]{.45\linewidth}
    \centering
    \includegraphics[width=.95\linewidth]{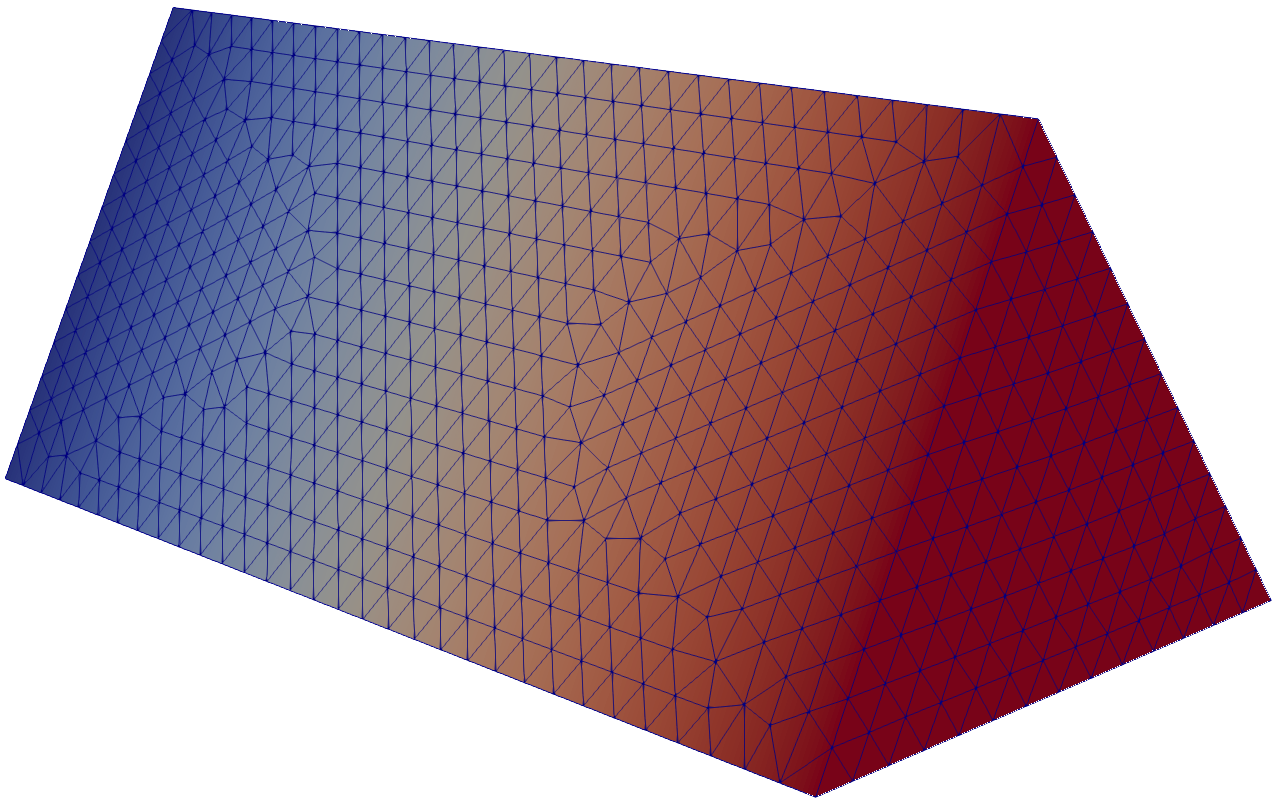}
    \caption{A refined skeleton mesh with 3368 elements.}
  \end{subfigure}
  \hfill
  \begin{subfigure}[b]{.45\linewidth}
    \centering
    \includegraphics[width=.95\linewidth]{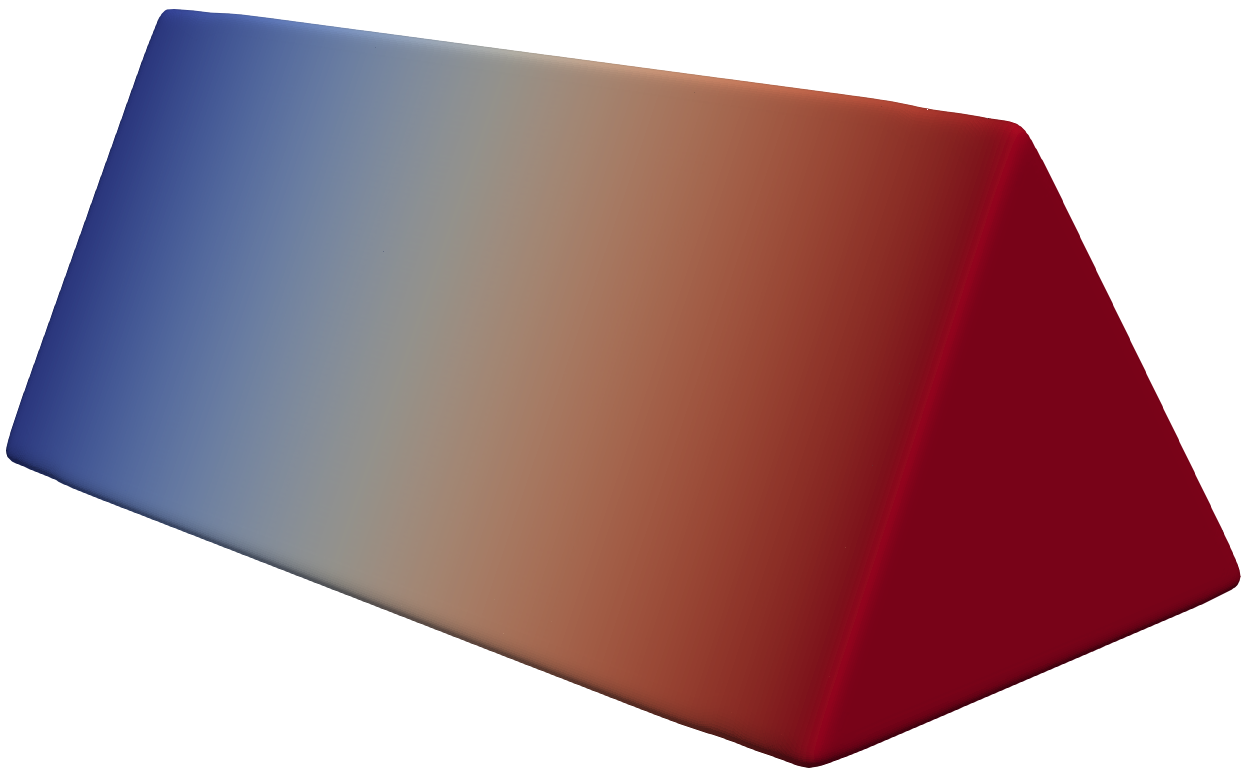}
    \caption{A smoothing of the refined mesh with~$\lambda=5$.}
  \end{subfigure}
  \caption{The main mode of failure of the algorithm of this work,
    so-called \emph{geometric ringing}. Qualitatively, the skeleton
    mesh is not a good approximation of an underlying smooth surface
    in the top figures. With some refinement, as in the bottom
    figures, the results are qualitatively more pleasing.}
  \label{fig:prism}
\end{figure}

When the mesh is globally coarse, it can be shown that the level set
function~$\Phi$  develops oscillations. This is most easily seen in
2D by examining the level sets (curves) of the convolution of a
Gaussian with the indicator function of the unit square.
We refer to this behavior as \emph{geometric ringing}, inspired by
similar behavior that arises in signal processing applications
(i.e. Fourier aliasing and filtering).
A simple three-dimensional example is shown in Figure~\ref{fig:prism}.
Notice also, in Figure~\ref{fig:prism_ringing},
that the volume/dimensions of the  object shrink noticeably
when a very coarse mesh is used as a skeleton surface. The mesh was
sufficiently coarse in this example that the algorithm failed to find
a level set during the Newton iterations with~$\lambda =
2.5$. Setting~$\lambda=10$ resulted in convergence to a level set, and
preserved more local feature of the geometry.

\section{Conclusions}
\label{sec:conclusions}

In this work, we have presented a multiscale algorithm which is able to 
construct arbitrarily high-order piecewise smooth surfaces using only
a skeleton mesh. In practice, this usually consists of a water-tight
flat triangulation as input. The algorithm uses a target-dependent
smoothing kernel, which is designed to be
commensurate with the local mesh size. The resulting smooth surface is constructed as
an atlas of charts from the skeleton triangles to a well-defined level set.
The scheme is accelerated via reformulation as a boundary integral and
the use of a fast multipole method; several numerical experiments illustrate
the performance of the approach.

It is worth reminding the reader of the principal mode of failure of the
algorithm: significant geometric artifacts will be generated when the
initial skeleton mesh does not accurately represent an underlying
smooth surface. That is to say, when the skeleton triangles are
not scaled to be commensurate with rapid changes in geometric features 
(such as the normal direction), we observe what we referred to as \emph{geometric
  ringing}, in which the smooth surface develops oscillations. This
behavior can be ameliorated by refining the skeleton
mesh to the desired length scale of smoothing.

Going forward, there are two important scenarios that we have not
addressed, but that are the focus of current work: (1) extending the algorithm
to skeleton meshes that are \emph{not} water-tight, and noisy
(e.g. those obtained from a 3D scanning procedure), and (2) permitting
the inclusion of singular features (such as corners and edges) which are to 
be left intact as part of the smoothing process. The algorithm, as it is
currently implemented, smooths away geometric singularities at the length
scale of the triangles impinging on the relevant feature. We are also 
investigating the use of our surface smoothing framework for problems involving
surface motion or surface remeshing.

\section*{Acknowledgments}
 We gratefully acknowledge the support of the NVIDIA Corporation 
 with the donation of a Quadro P6000, used for some of the 
 visualizations presented in this research.

\bibliographystyle{abbrv}
\bibliography{master}

\begin{thebibliography}{10}

\bibitem{ahrens2005paraview}
J.~Ahrens, B.~Geveci, and C.~Law.
\newblock Paraview: An end-user tool for large data visualization.
\newblock {\em The visualization handbook}, 717, 2005.

\bibitem{blinn1982}
J.~F. Blinn.
\newblock {A Generalization of Algebraic Surface Drawing}.
\newblock {\em ACM Trans. Graphics}, 1:235--256, 1982.

\bibitem{bloomenthal1991}
J.~Bloomenthal and K.~Shoemake.
\newblock Convolution surfaces.
\newblock {\em SIGGRAPH Comput. Graph.}, 25(4):251--256, July 1991.

\bibitem{borm2003introduction}
S.~B{\"o}rm, L.~Grasedyck, and W.~Hackbusch.
\newblock Introduction to hierarchical matrices with applications.
\newblock {\em Engineering Analysis with Boundary Elements}, 27(5):405--422,
  2003.

\bibitem{bremer-2015}
J.~Bremer, A.~Gillman, and P.-G. Martinsson.
\newblock A high-order accelerated direct solver for integral equations on
  curved surfaces.
\newblock {\em BIT Num. Math.}, 55:367--397, 2015.

\bibitem{bruno2001fast}
O.~P. Bruno and L.~A. Kunyansky.
\newblock A fast, high-order algorithm for the solution of surface scattering
  problems: {B}asic implementation, tests, and applications.
\newblock {\em J. Comput. Phys.}, 169(1):80--110, 2001.

\bibitem{wideband3d}
H.~Cheng, W.~Y. Crutchfield, Z.~Gimbutas, L.~Greengard, J.~F. Ethridge,
  J.~Huang, V.~Rokhlin, N.~Yarvin, and J.~Zhao.
\newblock A wideband fast multipole method for the {H}elmholtz equation in
  three dimensions.
\newblock {\em J. Comput. Phys.}, 216:300--325, 2006.

\bibitem{fmm2}
H.~Cheng, L.~Greengard, and V.~Rokhlin.
\newblock A fast adaptive multipole algorithm in three dimensions.
\newblock {\em J. Comput. Phys.}, 155(2):468--498, 1999.

\bibitem{coifman1993fast}
R.~Coifman, V.~Rokhlin, and S.~Wandzura.
\newblock The fast multipole method for the wave equation: A pedestrian
  prescription.
\newblock {\em IEEE Antennas Propag. Mag.}, 35(3):7--12, 1993.

\bibitem{cottrell}
J.~A. Cottrell, T.~J.~R. Hughes, and Y.~Bazilevs.
\newblock {\em Isogeometric Analysis: Toward Integration of CAD and FEA}.
\newblock Wiley, West Sussex, UK, 2009.

\bibitem{dapogny2014remesh}
C.~Dapogny, C.~Dobrzynski, and P.~Frey.
\newblock Three-dimensional adaptive domain remeshing, implicit domain meshing,
  and applications to free and moving boundary problems.
\newblock {\em J. Comput. Phys.}, 262:358--378, 2914.

\bibitem{darve2004fast}
E.~Darve and P.~Hav{\'e}.
\newblock {A fast multipole method for Maxwell equations stable at all
  frequencies}.
\newblock {\em Philosophical Transactions of the Royal Society of London A:
  Mathematical, Physical and Engineering Sciences}, 362(1816):603--628, 2004.

\bibitem{derose1998subdivision}
T.~DeRose, M.~Kass, and T.~Truong.
\newblock Subdivision surfaces in character animation.
\newblock In {\em Proceedings of the 25th annual conference on Computer
  graphics and interactive techniques}, pages 85--94. ACM, 1998.

\bibitem{epstein_2016}
C.~L. Epstein and M.~O'Neil.
\newblock Smoothed corners and scattered waves.
\newblock {\em SIAM J. Sci. Comput.}, 38:A2665--A2698, 2016.

\bibitem{fleishman2005}
S.~Fleishman, D.~Cohen-Or, and C.~T. Silva.
\newblock Robust moving least-squares fitting with sharp features.
\newblock {\em ACM Trans. Graph.}, 24(3):544--552, July 2005.

\bibitem{darvebb}
W.~Fong and E.~Darve.
\newblock The black-box fast multipole method.
\newblock {\em J. Comput. Phys.}, 228(23):8712--8725, 2009.

\bibitem{friedrichs}
K.~O. Friedrichs.
\newblock On the differentiability of the solutions of linear elliptic
  differential equations.
\newblock {\em Commun. Pure Appl. Math.}, 6:299--325, 1953.

\bibitem{gmsh}
C.~Geuzaine and J.-F. Remacle.
\newblock Gmsh: {A} 3-{D} finite element mesh generator with built-in pre- and
  post-processing facilities.
\newblock {\em Int. J. Num. Methods Engrg.}, 79:1309--1331, 2009.

\bibitem{ggm1993}
A.~Greenbaum, L.~Greengard, and G.~B. McFadden.
\newblock Laplace's equation and the {D}irichlet-{N}eumann map in multiply
  connected domains.
\newblock {\em J. Comput. Phys.}, 105(2):267--278, 1993.

\bibitem{greengard-1987}
L.~Greengard and V.~Rokhlin.
\newblock {A Fast Algorithm for Particle Simulations}.
\newblock {\em J. Comput. Phys.}, 73:325--348, 1987.

\bibitem{greengard-1997}
L.~Greengard and V.~Rokhlin.
\newblock {A new version of the Fast Multipole Method for the Laplace equation
  in three dimensions}.
\newblock {\em Acta Numerica}, 6:229--269, 1997.

\bibitem{hao_2014}
S.~Hao, A.~H. Barnett, P.-G. Martinsson, and P.~Young.
\newblock High-order accurate {N}ystr\"om discretization of integral equations
  with weakly singular kernels on smooth curves in the plane.
\newblock {\em Adv. Comput. Math.}, 40:245--272, 2014.

\bibitem{helsing}
J.~Helsing and R.~Ojala.
\newblock {C}orner singularities for elliptic problems: {I}ntegral equations,
  graded meshes, quadrature, and compressed inverse preconditioning.
\newblock {\em J. Comput. Phys.}, 227(20):8820--8840, 2008.

\bibitem{zorin2018}
Y.~Hu, Q.~Zhou, X.~Gao, A.~Jacobson, D.~Zorin, and D.~Panozzo.
\newblock Tetrahedral meshing in the wild.
\newblock {\em ACM Trans. Graph.}, 37:Article 60, 2018.

\bibitem{hughes2005isogeometric}
T.~J. Hughes, J.~A. Cottrell, and Y.~Bazilevs.
\newblock {Isogeometric analysis: CAD, finite elements, NURBS, exact geometry
  and mesh refinement}.
\newblock {\em Computer methods in applied mechanics and engineering},
  194(39-41):4135--4195, 2005.

\bibitem{gid}
{International Center for Numerical Methods in Engineering (CIMNE) }.
\newblock {GiD: The personal pre and post processor}.
\newblock \url{www.gidhome.com}, 2017.

\bibitem{koornwinder_1975}
T.~Koornwinder.
\newblock Two-variable analogues of the classical orthogonal polynomials.
\newblock In {\em Theory and application of special functions (Proc. Advanced
  Sem., Math. Res. Center, Univ. Wisconsin, Madison, Wis., 1975)}, pages
  435--495. Academic Press New York, 1975.

\bibitem{martinsson-2005}
P.-G. Martinsson and V.~Rokhlin.
\newblock {A fast direct solver for boundary integral equations in two
  dimensions}.
\newblock {\em J. Comput. Phys}, 205:1--23, 2005.

\bibitem{oneil2018surface}
M.~O'Neil.
\newblock {Second-kind integral equations for the Laplace-Beltrami problem on
  surfaces in three dimensions}.
\newblock {\em Adv. Comput. Math.}, 44(5):1385--1409, 2018.

\bibitem{phillips1997}
J.~Phillips and J.~White.
\newblock A precorrected-{FFT} method for electrostatic analysis of complicated
  3-{D} structures.
\newblock {\em IEEE Trans. Computer-Aided Design}, 16(10):1059--1072, 1997.

\bibitem{rachh_2016}
M.~Rachh, A.~Kl\"ockner, and M.~O'Neil.
\newblock Fast algorithms for {Q}uadrature by {E}xpansion {I}: {G}lobally valid
  expansions.
\newblock {\em J. Comput. Phys.}, 345:706--731, 2017.

\bibitem{freecad}
J.~Riegel, W.~Mayer, and Y.~van Havre.
\newblock {FreeCAD}, v. 0.18.1, 2019.
\newblock \url{http://www.freecadweb.org/}.

\bibitem{sherstyuk1999design}
A.~Sherstyuk.
\newblock Interactive shape design with convolution surfaces.
\newblock In {\em Proceedings Shape Modeling International '99. International
  Conference on Shape Modeling and Applications}, pages 56--65, March 1999.

\bibitem{sherstyuk1999}
A.~Sherstyuk.
\newblock Kernel functions in convolution surfaces: a comparative analysis.
\newblock {\em Visual Comput.}, 15:171--182, 1999.

\bibitem{Siegel2018}
M.~Siegel and A.-K. Tornberg.
\newblock {A local target specific quadrature by expansion method for
  evaluation of layer potentials in 3D}.
\newblock {\em J. Comput. Phys.}, 364:365--392, 2018.

\bibitem{simpson2014acoustic}
R.~N. Simpson, M.~A. Scott, M.~Taus, D.~C. Thomas, and H.~Lian.
\newblock Acoustic isogeometric boundary element analysis.
\newblock {\em Computer Methods in Applied Mechanics and Engineering},
  269:265--290, 2014.

\bibitem{song1997multilevel}
J.~Song, C.-C. Lu, and W.~C. Chew.
\newblock Multilevel fast multipole algorithm for electromagnetic scattering by
  large complex objects.
\newblock {\em IEEE Trans. Antennas Propag.}, 45(10):1488--1493, 1997.

\bibitem{strain91}
J.~Strain.
\newblock The fast {G}auss transform with variable scales.
\newblock {\em SIAM Journal on Scientific and Statistical Computing},
  12(5):1131--1139, 1991.

\bibitem{vioreanu_2014}
B.~Vioreanu and V.~Rokhlin.
\newblock {Spectra of Multiplication Operators as a Numerical Tool}.
\newblock {\em SIAM J. Sci. Comput.}, 36:A267--A288, 2014.

\bibitem{wala2018qbx}
M.~Wala and A.~Kl\"ockner.
\newblock {A Fast Algorithm for Quadrature by Expansion in Three Dimensions}.
\newblock {\em arXiv [math.NA]}, 1805.06106, 2018.

\bibitem{wang2018}
J.~Wang and L.~Greengard.
\newblock An adaptive fast {G}auss transform in two dimensions.
\newblock {\em SIAM J. Sci. Comput.}, 40:A1274--A1300, 2018.

\bibitem{ying}
L.~Ying, G.~Biros, and D.~Zorin.
\newblock A high-order 3{D} boundary integral equation solver for elliptic
  {PDE}s in smooth domains.
\newblock {\em J. Comput. Phys.}, 219(1):247--275, 2006.

\bibitem{ying2004simple}
L.~Ying and D.~Zorin.
\newblock A simple manifold-based construction of surfaces of arbitrary
  smoothness.
\newblock {\em ACM Trans. Graphics}, 23(3):271--275, 2004.

\bibitem{zorin1996interpolating}
D.~Zorin, P.~Schr{\"o}der, and W.~Sweldens.
\newblock Interpolating subdivision for meshes with arbitrary topology.
\newblock In {\em Proceedings of the 23rd annual conference on Computer
  graphics and interactive techniques}, pages 189--192. ACM, 1996.

\end{thebibliography}

\end{document}